\documentclass[11pt]{article}
\usepackage{subcaption}
\usepackage{amsmath}
\usepackage{amssymb}
\usepackage{latexsym}
\usepackage{amsmath, amsfonts,amssymb, amsthm, euscript,makeidx, color, mathrsfs}
\usepackage{xcolor}
\usepackage{indentfirst}
 \usepackage{algorithm}
  \usepackage{booktabs}
 \usepackage[noend]{algpseudocode}
\usepackage{graphicx}
\usepackage{epstopdf}
\usepackage{float}
\usepackage{subfig}
\usepackage{multirow}
\newtheorem{assumption}{Assumption}
\def\qed{ \ \vrule width.2cm height.2cm depth0cm\smallskip}

\newcommand{\ba}{\begin{array}}
\newcommand{\ea}{\end{array}}
\newcommand{\be}{\begin{equation}}
\newcommand{\ee}{\end{equation}}
\newcommand{\bea}{\begin{eqnarray}}
\newcommand{\eea}{\end{eqnarray}}
\newcommand{\beaa}{\begin{eqnarray*}}
\newcommand{\eeaa}{\end{eqnarray*}}

\def\neg{\negthinspace}

\def\a{\alpha}
\def\g{\gamma}
\def\d{\delta}
\def\e{\varepsilon}

\def\l{\lambda}
\def\m{\mu}

\def\si{\sigma}
\def\t{\tau}

\def\th{\theta}
\def\o{\omega}
\def\h{\widehat}

\def\D{\Delta}

\def\O{\Omega}

\def\Th{\Theta}

\def\D{\Delta}
\def\Th{\Theta}

\def\O{\Omega}
\def\cA{{\cal A}}

\def\cF{{\cal F}}

\def\cH{{\cal H}}

\def\hC{\mathbb{C}}
\def\hD{\mathbb{D}}
\def\hE{\mathbb{E}}
\def\hF{\mathbb{F}}

\def\hH{\mathbb{H}}

\def\hL{\mathbb{L}}

\def\hN{\mathbb{N}}

\def\hP{\mathbb{P}}
\def\hQ{\mathbb{Q}}
\def\hR{\mathbb{R}}

\def\hX{\mathbb{X}}

\def\sA{\mathscr{A}}
\def\sB{\mathscr{B}}

\def\sP{\mathscr{P}}

\def\sU{\mathscr{U}}

\def\no{\noindent}

\def\ss{\smallskip}
\def\ms{\medskip}
\def\bs{\bigskip}
\def\q{\quad}
\def\qq{\qquad}

\def\pa{\partial}
\def\cd{\cdot}
\def\cds{\cdots}

\def\qed{ \hfill \vrule width.25cm height.25cm depth0cm\smallskip}

\newcommand{\basa}{\begin{assumption}}
\newcommand{\easa}{\end{assumption}}

\newcommand{\ol}{\overline}

\newcommand{\bas}{\begin{assum}}
\newcommand{\eas}{\end{assum}}

\def\limsup{\mathop{\overline{\rm lim}}}
\def\liminf{\mathop{\underline{\rm lim}}}

\def\pa{\partial}
\def\h{\widehat}

 \def\cd{\cdot}
\def\cds{\cdots}

\def\ol{\overline}
\def\ul{\underline}

\def\dis{\displaystyle}

\def\1{{\bf 1}}

\def\ppi{{\pmb\pi}}

\def\:{\!:\!}
\def\reff#1{{\rm(\ref{#1})}}

at 9pt

\def\prehp(#1,#2){\ensuremath{  #1 \cdot #2 }}

\begin{document}

\newtheorem{thm}{Theorem}[section]
\newtheorem{lem}[thm]{Lemma}
\newtheorem{cor}[thm]{Corollary}
\newtheorem{prop}[thm]{Proposition}
\newtheorem{rem}[thm]{Remark}
\newtheorem{eg}[thm]{Example}
\newtheorem{defn}[thm]{Definition}
\newtheorem{assum}[thm]{Assumption}

\renewcommand {\theequation}{\arabic{section}.\arabic{equation}}
\def\thesection{\arabic{section}}

%\begin{document}

\title{\bf Reinforcement Learning for Optimal Dividend Problem under Diffusion Model }

\author{Lihua Bai\thanks{\noindent School of Mathematics,
Nankai University, Tianjin, 300071, China. Email:
lbai@nankai.edu.cn. This author is supported in part by Chinese NSF
grants \#12272274 and \#12571514.}, ~Thejani Gamage\thanks{\noindent
Department of Mathematics and Statistics, 
University of Massachusetts Amherst, Amherst, MA 01003. Email: tgamage@umass.edu. }, ~Jin Ma\thanks{\noindent
 Department of Mathematics,
University of Southern California, Los Angeles, CA 90089.
%; School of Statistics and Mathematics, Shanghai Lixin University of Accounting and Finance, Shanghai, China.
Email: jinma@usc.edu. This author is supported in part by NSF grant \#2205972. }, ~Gaozhan Wang \thanks{\noindent Department of Mathematics,
University of Southern California, Los Angeles, CA 90089. Email: gaozhanw@usc.edu.}}
%The correct dates will be entered by the editor.
\date{}
\maketitle

\begin{abstract}
In this paper, we study the optimal dividend problem under the continuous time diffusion model with the bounded dividend rate from the 
{\it Reinforcement Learning (RL)} perspective. Unlike the standard literature, our main focus will be on  numerical algorithms that allow part or all of the  system parameters to be unspecified
% issue of, that is, when the  {\it parameter uncertainty} 
so that the optimal control cannot be explicitly determined. 
Following the RL literature, we introduce the  {\it entropy-regularized exploratory control problem}, which randomizes the control actions and balances the levels of exploitation and exploration, and carry out a theoretical analysis of the associated 
%its value function and related HJB equation.  We will then use a 
 {\it Policy Improvement} (PI) and  {\it Policy Evaluation} (PE) devices
% {\it Policy Improvement} (PI) procedure as well as the  {\it Policy Evaluation} (PE) devices 
% via the 
% %Temporal Difference (TD)-based algorithm  and 
% Martingale Loss (ML) function introduced in \cite{Zhou}
and the corresponding sequence of the approximating optimal strategies. Specifically, our algorithm will be based on two independent neural networks that approximate the value function and its derivative simultaneously. Such an algorithm, to the best of our knowledge, is new in the context of the optimal dividend problems, and can be effective even for the situation when the premium and/or interest rate is state-dependent. 
% hence beyond the reach of the standard statistical methods. 
 Multiple numerical experiments are presented to empirically demonstrate the effectiveness of our RL algorithm, especially the benefit of invoking two independent neural networks.
% Some numerical experiments are presented as well. }
%using parametrization families for the cost functional, to illustrate the effectiveness of the approximation schemes and to discuss the factors that affect the accuracy and effectiveness of the policy evaluation algorithms.  
\end{abstract}

\vfill \bs

\no


\bs

\no{\it 2020 AMS Mathematics subject classification:} {93E20,35; 91G50; 93B47.}

% 60H05: stochastic integrals
% 60H10: stochastic ordinary differential equations
% 60G44: martingales with continuous parameter
% 28B20: set-valued set functions and measures; integration of set-valued functions; measurable selections
% 47H04: set-valued operators

\eject

\section{Introduction}
The problem of maximizing the cumulative discounted dividend payment can be traced back to the work of de Finetti \cite{de}.
In the continuous time diffusion case, Asmussen and Taksar \cite{asmussen} solved the optimal dividend problem for the first time, in both constrained rate % They considered two situations. In the first situation, 
(i.e., the dividend rate is restricted to a finite interval $[0, a]$), and in the general unbounded rate case. The optimal control was shown to be of a threshold type in the former and a barrier type in the latter, respectively. Other works on the
optimal dividend problem for controlled diffusion models include, e.g.,  \cite{PaGj, RaSh, Tak-Zhou}, and many variations can be found in, e.g.,  \cite{BP, CGY,fh2, DTE,  TA, ZY, fh4}, and the references cited therein.
It should be noted that in all the works above, the model was assumed to be ``certain" in the sense that all the parameters (e.g., drift and volatility) of the diffusion were specified, and the optimal control was a deterministic function of the given parameters. But in reality, the distribution of the claim process is usually unknown; thus, neither are the drift and volatility coefficients. 
%, as they depend on either the mean or variance of the claims. 
{\color{black}Thus, finding an efficient way to determine the optimal strategy allowing a certain degree of  ``model uncertainty" or ``parameter uncertainty", i.e., without specific knowledge of all the system parameters, is a longstanding challenge. }

The main motivation of this paper is to study the optimal dividend problems in a data-driven approach, and in particular,  when the model parameters are not specified, so the optimal control cannot be explicitly determined. More specifically, we shall approach this problem using the now popular {\it Reinforcement Learning} (RL) method, 
%especially a model-free  approach 
in which the agent observes multiple simulated trajectories of the state process to approximate the optimal strategy (see, e.g.,  \cite[Chapter 4]{HUM}). We should note that using reinforcement learning to solve discrete Markov decision problems has been well studied. 
% but the extension of these concepts to the continuous time and space setting is still fairly new. 
Roughly speaking, in RL the learning agent uses a sequence of \textit{trial and error} to determine the optimal control and optimal value function. Such a process has been characterized by a mixture of  \textit{exploration} and \textit{exploitation}, 
%and a critical point is to balance the  levels of each aspect, commonly by 
balanced through the ``randomized" control action and the additional Shannon's entropy weighted by a temperature parameter in the cost functional. Reinforcement Learning in continuous time and space has been studied by many authors (cf. e.g., \cite{JOZ, zhou,jia2022policy,  zhou1,  zhou2}) and the resulting control problem involving Shannon's entropy was often referred to as the {\it entropy-regularized exploratory optimal control problem}. The specific formulation of our particular problem   
% that allows us to use the model-free approach 
will be discussed in \S2.

Two key elements in the well-known {\it Policy Iteration Algorithms} (PIA) for solving an entropy-regularized exploratory optimal control problem are known as the {\it Policy Evaluation}  (PE) and {\it Policy Improvement}  (PI)  procedures. The former evaluates a given policy, and the latter produces a ``better" one. We should note that the idea of PIA, as well as its convergence analysis, is not new in the numerical optimal control literature (see, e.g.,  \cite{PIA1, PIA2, PIA3, PIA4}). The main challenge in exploratory control problems, especially for continuous time diffusion models, is often the technical complications caused by the involvement of entropy regularization. There have been strong efforts to prove convergence, even the rate of convergence, in both finite and infinite horizon cases using policy iteration or other RL approaches 
% when the  structure of the optimal exploratory control is known
 (cf. e.g., \cite{HWZcvg, jia2022policy,MWZ1,sethi2025entropy,   TWZ}). Although our scheme is quite close to the existing infinite horizon PIA, there are actually a few differences.

First, our cost functional involves the ruin time, which complicates the problem significantly. In particular, the approaches in the aforementioned works do not seem to work easily here, due to the irregularity of the ruin time with respect to the initial states, as it is not differentiable in general. 
Our theoretical analysis takes advantage of the fact that in this paper the state process is one-dimensional, taking nonnegative values, so that some stability arguments for the first-order nonlinear systems can be applied to conclude that the exploratory HJB equation has a concave, bounded classical solution, whence a viscosity solution that coincides with the value function of the optimal dividend problem. Second, to establish a robust and stable numerical RL algorithm, we follow a recently developed neural network approach, which is particularly effective in the situation when the model parameters 
%(especially the drift coefficients of the diffusion models) 
are not known (cf. \cite{MWZZ26}). {\color{black} A key element in such an algorithm is to utilize an independent neural network via the so-called Bismut-Elworthy-Li(BEL) representation formula (see, e.g., \cite{Bismut, EL, MaZhang1})
%, which utilizes two independent neural networks that 
to approximate the derivative of the value function directly, without invoking the regularity issue of the ruin time. 
%With the help of these regularity results, we prove the convergence of PIA with a much simpler argument.
%As we argued in \S7 of the revised version, we can somewhat embed it into an infinite horizon problem but defined on $\hR_+$. The technical novelty of our final RL algorithm is the use of, which allows us to establish an independent neural network 
To the best of our knowledge, such an approach to dealing with the state-dependent random horizon (ruin time) is novel. 
Furthermore, it is worth noting that, unlike the traditional optimal dividend problem, the true solution of the entropy-regularized HJB in this paper cannot be explicitly solved due to the transcendental nature of the equation. Therefore, the convergence of our algorithm will be assessed through the performance of both PE (vanishing of martingale loss) and PI (the improvement of the value function). But on the other hand, although our problem is ``stopped" at a ruin time, it is essentially an infinite horizon problem by nature, which somewhat prevents the use of the so-called “batch” learning method that depends on multiple (whole) trajectories. In our algorithm, we consider  
%Nevertheless, we can still study 
the temporally “truncated” problem so that the batch learning method can be applied.   It turns out that the truncated trajectories can produce the desired numerical results with negligible error when the discounting factor is relatively large, as 
%for a suitable truncation time $T$, because 
the truncation error by nature has an exponential decay.

Finally, in the numerical section (\S6) we carry out experiments to empirically validate the PIA
scheme, both in the constant coefficient case and with a state-dependent drift to test the general applicability of our algorithm. 
% by stracking the decay of the policy evaluation loss and the improvement of the value function.  namely with constant coefficients and with a state-dependent drift coefficient,  
In both cases, we observe satisfactory
empirical convergence of the proposed PIA scheme.  We also discuss the financial implications of the approximated optimal strategy. 
%, which cannot be easily obtained theoretically from the entropy-regularized HJB equation itself, due to its transcendental nature. 
%Furthermore, we 
Our algorithm features two independent   neural networks  that approximates the value function and its derivative simultaneously   via  the BEL representation, and we compare our numerical results with those that approximate the derivatives
%and explain the theoretical rationale for approximating
%the derivative of the value function by an independent neural
%network via, rather than obtaining it
solely through {\it automatic differentiation} (Autodiff) of the
learned value function. It is worth noting that the 
%The key distinction is that 
convergence   in $\hL^2$ does not, in general, imply the
convergence in $\hH^1$, thus an accurate approximation of the value function through martingale loss
alone does not guarantee accurate  approximation of its derivative (see also \cite{MWZZ26}). This motivates training the derivative network directly against
its own representation-based target. We present both theoretical discussion and various numerical evidence showing that our algorithm does out-performed the one using Autodiff, in terms of both error reduction and sensitivity  to the truncation horizon and neural network structure.  

%Example~\ref{EX3} complements, showing
%improved derivative accuracy over Autodiff in the fixed-policy
%benchmark considered. 

}

{\color{black}
The rest of the paper is organized as follows. In \S2  we give the preliminary description of the problem and all the necessary notations, definitions, and assumptions. In \S3 we study the value function and its regularity, and prove that it is a concave, bounded classical solution to the exploratory HJB equation. In \S4 we study the policy update procedure, where we introduce our Policy Improvement scheme, and prove its convergence.
In \S5 we discuss the Martingale Loss-based approach for Policy Evaluation.
%, that is, the methods for approximating the cost functional for a given policy, using a 
 % and   (online) CTD($\g$) methods, respectively.
In \S6 we present our RL algorithm combining PE and PI (the PIA) and provide the results of numerical experiments. 
%a probabilistic representation based approach for the derivative estimation propose parametrization families for PE and present numerical experiments using the proposed PI and PE methods. %\section{Prelimilaries}
%\label{sec:Prelimilaries}
}

\section{Preliminaries and Problem Formulation}
\setcounter{equation}{0}

%\subsection{\bf Exploratory optimal dividend problem}
Throughout this paper we consider a filtered probability space $(\O, \cF, \{\cF_t\}_{t\geq0}, \hP)$ on which is defined a standard Brownian motion  $\{W_t,t\geq 0\}$. We assume that the filtration $\hF:=\{\cF_t\}=\{\cF^W_t\}$, with the usual augmentation so that it satisfies the usual conditions. For any metric space $\hX$ with topological Borel sets $\sB(\hX)$, we denote $\hL^0(\hX)$ to be all $\sB(\hX)$-measurable functions, and $\hL^p(\hX)$, $p\ge 1$, to be the space of $p$-th integrable functions. The spaces  $\hL^0_\hF([0,T];\hR)$
and $\hL^p_\hF([0,T];\hR)$, $p\ge1$, etc., are defined in the usual ways. Furthermore, for a given domain $\hD\subset \hR$, we denote $\hC^k(\hD)$ to be the space of all $k$-th order continuously differentiable functions on $\hD$, and $\hC(\hD)=\hC^0(\hD)$. In particular, for $\hR_+:=[0,\infty)$, we denote $\hC^k_b(\hR_+)$ to be the space of all bounded and $k$-th continuously differentiable functions on $\hR_+$ with all derivatives being bounded. 

Consider the simplest  diffusion approximation of a Cramer-Lundberg model with dividend:
\bea
\label{X0}
dX_t=(\mu-\a_t)dt+\sigma dW_t,\quad t>0,\quad X_0=x\in\mathbb{R},
\eea
where $x$ is the initial state,  $\mu$ and $\sigma$ are constants determined by the premium rate and the claim frequency and size (cf., e.g., \cite{asmussen}), and $\a_t$ is the dividend rate at time  $t\geq0$.  We denote $X=X^\a$ if necessary, and say that 
$\alpha=\{\a_t,t\geq0\}$ is {\it admissible} if it is  $\hF$-adapted and takes values in 
%a controlled state process $\{X_t,t\geq0\}$ is assumed to evolve like a Brownian motion .  A
a given  ``action space" $[0,a]$. Furthermore, let us define the {\it ruin time} to be $\tau^\alpha_x:=\inf\{t>0:X^\a_t<0\}$. Clearly, 
% classical controlled diffusion model with restricted dividend payment problem has the following state process:  If 
$X^\a_{\t^\a}=0$, and the problem is considered ``ruined" as no dividend will be paid after $\t^\a$.  Our aim is to maximize the expected total discounted dividend given the initial condition $X^\a_0=x\in\mathbb{R}$:
% and to obtain the value function:
\bea
\label{value}
 \qq V(x):=\sup_{\alpha\in \sU[0,a]}\hE_x\Big[\int_{0}^{\tau_x^\a} e^{-ct} \a_tdt\Big]:=\sup_{\alpha\in \sU[0,a]}\hE\Big[\int_{0}^{\tau^\a_x} e^{-ct} \a_tdt\Big|x_0=x\Big],
\eea
where $c>0$ is the discount rate, and $\sU[0,a]$ is the set of {\it admissible} dividend rates taking values in $[0,a]$. The problem (\ref{X0})-(\ref{value}) is often referred to as the {\it classical optimal dividend} problem with the {\it restricted dividend rate} in a given interval $[0, a]$. 
%  and $a_t$ is called reward in the RL formulation. 

It is well-understood that when the parameters $\mu$ and $\sigma$ are known, then the optimal control is of the ``feedback" form:  $\a_t^*=\boldsymbol{a}^*(X_t^*)$, where $X_t^*$ is the corresponding state process   and $\boldsymbol{a}^*(\cd)$ is a deterministic function taking values in $[0,a]$, often in the form of  a {\it threshold control} (see, e.g.,  \cite{asmussen}). However, in practice 
%a serious problem is that 
the exact form of $\boldsymbol{a}^*(\cd)$ is not implementable since the model parameters are usually not known; thus the ``parameter insensitive"  method through {\it Reinforcement Learning} (RL) becomes a much more desirable alternative,  which we now elaborate. 

In the RL formulation, the agent follows a process of exploration and exploitation via a sequence of trial-and-error evaluation.  A key element is to randomize the control action as a probability distribution over $[0, a]$, 
% and assume that all the trials are independent and identically distributed, very much 
similar to the notion of {\it relaxed control} in control theory, and the classical control is considered as a special point-mass (or Dirac $\d$-measure)  case.
To make the idea more accurate mathematically, let us
denote $\sB([0,a])$ to be the Borel field on $[0,a]$, and $\sP([0,a])$ to be the space of all probability measure on $([0,a], \sB([0,a]))$, endowed with, say, the Wasserstein metric. 
A ``relaxed control" is  a randomized policy defined as a measure-valued progressively measurable process
$(t, \o)\mapsto { \pi}(\cd; t, \o)\in \sP([0,a])$. Assuming that  ${ \pi}(\cd;t, \o)$ has a density, denoted by $\pi_t(\cd, \o)\in\hL^1_+([0,a])\subset \hL^1([0,a])$, $(t, \o)\in[0,T]\times\O$, then we can write
$${  \pi}(A; t, \o)=\int_A \pi_t(w, \o)dw, \qq A\in \sB([0,a]), \q (t, \o)\in[0,T]\times \O.
$$
In what follows, we shall often identify a relaxed control with its density process  $\pi=\{\pi_t,t\geq 0\}$.
Now, for $t\in[0,T]$, we define a probability measure on $([0,a]\times\O, \sB([0, a])\otimes\cF)$ as follows: for  $ A\in\sB([0,a])$ and
$B\in\cF$,
\bea
\label{Qt}
\hQ_t(A\times B):=\int_A\int_B\pi(dw; t, \o)\hP(d\o)=\int_A\int_B\pi_t( w, \o)dw\hP(d\o).
\eea
%Let us define
%    $ \tilde{\mathcal{F}}_t = \sB([0,a])\otimes\cF_t$, $t\ge 0$, and consider the 
We call a function $A^\pi: [0,T]\times[0,a]\times\O\mapsto [0, a]$ the ``canonical representation" of a relaxed control $\pi=\{\pi(\cd, t, \cd)\}_{t\ge0}$, if $A^\pi_t(w, \o)=w$. Then, for $t\ge 0$ we have
%, the expectation of the random variable $A^\pi_t$  under $\hQ_t$ would be
\bea
\label{EQtA}
\hE^{\hQ_t}[A^\pi_t]=\int_\O\int_0^a A^\pi_t(w, \o)\pi(dw;t, \o)\hP(d\o)=\hE^{\hP}\Big[\int_0^a w\pi_t(w)dw\Big].
\eea
%In this way, the classical control can be viewed as a sample of the process $A^\pi$.

We can now derive the exploratory dynamics of the state process $X$ along the lines of entropy-regularized relaxed stochastic control arguments (see, e.g.,  \cite{zhou1}). Roughly speaking, consider the discrete version of the dynamics (\ref{X0}): for small $\D t>0$,
\bea
\label{zengliang}
\Delta x_t :=x_{t+\Delta t} -x_t \approx (\mu-a_t )\Delta t+\sigma (W_{t+\Delta t} -W_t ),\qquad t\geq0.
\eea
Let $\{a_t^i\}_{i=1}^N$ and  $\{(x^i_t,W_t^i)\}_{i=1}^N$  be $N$ independent samples of $(a_t)$ under the distribution $\pi_t$, 
%The superscript $\pi$ represents the control is under distribution denoted by $\pi=\{\pi_t(w),t\geq 0\}$.
and  the corresponding  samples of  $(X_t^\pi,W_t)$, respectively.
%, under the  controls $\{a^i_t\}_{t\ge 0}$, $i=1,2,...,N$ (or equivalently, $A^\pi$). 
%Then the increment of each process $x^i$ in the time interval $[t,t+\Delta t]$, with $\Delta t$ small enough, can be expressed as
%
Then, the law of large numbers and \eqref{EQtA} imply that
%, \eqref{zengliang} and , for $N$ large enough, we see that
\bea
\label{yczl}
\qq \sum_{i=1}^{N}\frac{\Delta x_t^i}{N}\approx \sum_{i=1}^{N} (\m-a^i_t)\frac{\D t}{N}\approx \hE^{\hQ_t}[\m-A^\pi_t]\D t=\hE^\hP\Big[\mu\neg-\neg\int_{0}^{a}w\pi_t (w, \cd)dw\Big]\Delta t,
  \eea
as $N\to\infty$. This, together with the fact $\frac{1}{N}\sum_{i=1}^{N}(\Delta x_t^i)^2\approx\sigma^2\Delta t$, leads to the following form 
%The identities \eqref{zengliang} and \eqref{yczl} motivate us to propose the 
of the exploratory version of the state dynamics:
\begin{eqnarray}\label{xfc}
dX_t=\Big(\mu-\int_{0}^{a}w\pi_t(w, \cd)dw\Big)dt+\sigma dW_t, \q X_0=x,
\end{eqnarray}
where $\{\pi_t(w, \cd)\}$ is the (density of) relaxed control process, and we shall often denote $X=X^{\pi, 0,x}=X^{\pi, x}$ to specify its dependence on control $\pi$ and
the initial state $x$.

To formulate the {\it entropy-regularized} optimal dividend problem, we first give a heuristic argument. Similar to (\ref{yczl}), for $N$ large and $\D t$ small we should have
\begin{eqnarray*}
\frac{1}{N}\sum_{i=1}^{N}e^{-ct}a_t^i\1_{[t\leq \tau^i]}\Delta t
%\xrightarrow{a.s.} 
\approx\hE^{\hQ_t}\Big[e^{-ct}A_t^\pi \1_{[t\leq \tau^\pi_x]}\Delta t\Big]
= \hE^\hP\Big[\1_{[t\leq \tau^\pi_x]}e^{-ct}\int_{0}^{a}w\pi_t(w)dw\Delta t\Big].
\end{eqnarray*}
{\color{black}
Therefore, in light of \cite{zhou1},  we shall define the entropy-regularized cost functional of the optimal expected dividend control problem under the  relaxed control ${\pi}$ with ``temperature $\l$" as 
\bea
\label{Jxpi}
&& J(x, {\pi)}=\hE_x\Big[\int_0^{\tau_x^{{\pi}}}e^{-ct} r^ \pi_\l(t) dt \Big],
\eea
where $ {\tau^{\pi}_x} = \inf  \{t > 0: X_t^{\pi,x} < 0  \}$, and 
$$r^ \pi_\l(t):= \int_{0}^{a}(w-\lambda\ln\pi_t(w))\pi_t(w)dw =\int_0^a w\ln \pi_t(w)dw+\l \cH(\pi). 
$$
Here, $\cH(\pi):=-\int_{0}^{a} \ln\pi_t(w))\pi_t(w)dw$  is the so-called {\it Shannon's entropy} and $\lambda>0$ is the {\it temperature parameter} balancing the exploration and exploitation.}

%{\color{blue} I think we should keep $r$ here instead of $\cH$ since $\cH$ is commonly used for entropy and we have an entropy-regularized reward here. We use two notations throughout the paper; we need to fix the notation. Let me know what you think a good notation should be, and I can make the notation consistent throughout.}

We now define the set of {\it open-loop admissible controls} as follows.
\begin{defn}
\label{Uad}
A measurable (density) process ${ \pi}=\{\pi_t(\cd,\cd)\}_{t\ge0}\in \hL^0([0, \infty)\times[0,a]\times\O)$ is called an open-loop admissible relaxed control if

1. ${{\pi}_t}(\cd;\o)\in\hL^1([0,a])$, for $dt\otimes d\hP$-a.e. $(t, \o)\in[0,\infty)\times \O$;

\ss
2. for each $w\in [0,a] $, the process $(t, \o)\mapsto {{\pi}_t}(w, \o)$ is $\hF$-progressively measurable;

\ss
 3.~ $\hE_x\big[\int_{0}^{\tau^\pi_x}e^{-ct} |r_\l^\pi (t) |dt\big]<+\infty$.

We shall denote $\sA(x)$ as the set of open-loop admissible relaxed controls.
\qed
\end{defn}

Consequently, we define 
%instead of the value function (\ref{value}), we now have 
the ``exploratory" value function as
\begin{eqnarray}
\label{v}
&& V(x)= \sup_{{\pi}\in\sA (x)}\hE_x\Big\{\int_{0}^{\tau_x^\pi}e^{-ct}r_\l^\pi(t) dt\Big\},
\qq x\ge 0.
\end{eqnarray}

An important type of   $\pi\in \sA(x)$ is of the  
``feedback" nature, that is, $\pi_t(w, \o)=\pmb\pi(w, X^{\pmb\pi,x}_t(\o))$ 
%$(t, \o, w)\in[0,T]\times \O\times [0,a]$ 
for some  deterministic function $\pmb\pi$, where 
%:[0,a]\times\hR\mapsto \hR$. Note that in this case, the state dynamics 
$X^{\pmb\pi,x}$ satisfies the SDE:
\bea
\label{sde2}
dX_t=\Big(\mu-\int_{0}^{a}w {\pmb\pi} (w, X_t)dw\Big)dt+\sigma dW_t, \q t\ge 0; \q X_0=x.
\eea
\begin{defn}
\label{Cad}
A function $\pmb \pi\in\hL^0([0, a]\times\hR)$ is called a closed-loop admissible relaxed control if, for every $x>0$,

\ss
1.  The SDE $(\ref{sde2})$ admits a unique strong solution $X^{\pmb{\pi},x}$, 

\ss
2. The process  $\pi=\{{{\pi}_t}(\cd;\o):= {\pmb\pi}(\cd,X^{\pmb\pi,x}_t(\o)); (t,\o)\in[0, T]\times\O \} \in \sA(x)$.

We denote $\sA_{cl}\subset \sA(x)$  to be the set of closed-loop admissible relaxed controls. 
\qed
\end{defn}

The following properties of the value function are straightforward. 
\begin{prop}
\label{assl}
Assume  $a>1$. Then the value function $V$   satisfies the following properties:

\ss
(1) $V(x)\ge V(y)$, if $x\ge y>0$;

\ss

(2)  $0\leq V(x)\leq \frac{\lambda \ln a+a}{c}$, $x\in\hR_+$.

\end{prop}

{\it Proof.} (1) Let $x\ge y$, and $\pi\in\sA(y)$. Consider  $\hat \pi_t(w, \o):=\pi_t(w, \o)\1_{\{t<\t^\pi_y(\o)\}}
+\frac{e^{\frac{w}{\l}}} {\l (e^{\frac{a}{\l}}-1)}
\1_{\{t\ge \t^\pi_y(\o)\}}$, $(t,w, \o)\in[0, \infty)\times[0,a]\times\O$.
Then, it is readily seen that $J(x, \hat \pi)\ge J(y, \pi)$, for $a>1$. Thus $V(x)\ge J(x, \hat\pi)\ge V(y)$, proving (1), as $\pi\in\sA(y)$ is arbitrary. 
%The monotonicity of the value function $V$ is obvious by definition, as for each $\pmb\pi\in\sA(x)$, $\t^{\pmb\pi}_x\ge \t^{\pmb\pi}_y$ for  $x\ge y$, 

(2) By definition $\int_{0}^{a} w\pi_t(w)dw \leq a$ and $-\int_{0}^{a} \ln\pi_t(w))\pi_t(w)dw \leq  \ln a$ by the well-known Kullback-Leibler divergence property. Thus, $r_\l^\pi(t) \leq \l \ln a +a$, and then $V(x) \leq \frac{\l \ln a +a}{c}$. On the other hand, since $ \hat{\pmb \pi}(w,x) \equiv\frac{1}{a}$, $(x,w) \in \hR_+ \times [0,a]$, is admissible and $J(x,\hat{\pmb \pi}) \geq 0$ for $a>1$, the conclusion follows. 

\qed

 We remark that in optimal dividend control problems it is often assumed that the maximal dividend rate is greater than the average return rate (that is,  $a>2\mu$),  % has some interesting implications. For example, as a general purpose in the RL approach (see, e.g., \cite{rlshu1}), one often pays attention to the asymptotic behavior as $\l\to \infty$.
%But note that $e^{\frac{a}\l}=1+\frac{a}{\l}+o(\frac1\l)\sim 1+\frac{a}{\l}$, for  large $\l$, 
%the value function would be non-negative whenever $a>1$.
%, which we shall do in the sequel.
%ofencourage the state process to ruin as quickly as possible when we give non-positive reward for every time that from initial time to the end of state process. Then the proposition \ref{assl} motivates us to propose $a>1$. It is easy to accept as a fact when we change the unit of the state process.
%Furthermore, as  Combining,  
and that the average return of a surplus process $X$,
% of the insurance company to have a positive return that the return, 
including the safety loading, is higher than the interest rate $c$. These, together with Proposition \ref{assl}, 
lead to the following standing assumption that
%This, together with the analysis on the dividend rate, leads to the following assumption, which 
will be used throughout the paper. 

\begin{assum}
\label{jc}
(i) The maximal dividend rate $a$ satisfies  $a>\max\left\{1,2\mu\right\}$; and 

\ms
(ii) the average return $\mu$ satisfies $\mu >\max\left\{c,\sigma^2/2\right\}$.
\qed
\end{assum}

\section{The Value Function and Its Regularity}
\setcounter{equation}{0}

%\subsection{\bf Hamilton-Jacobi-Bellman Equation}

In this section we study the value function of the relaxed control problem (\ref{v}). We note that while most of the results are well-understood, some details still require justification, especially concerning the regularity, due particularly to the non-smoothness of the exit time $\t_x$.  

We begin by recalling the Bellman optimality principle (cf. e.g.,  \cite{YongZhou}): 
%, applying the stochastic version of Bellman's principle of optimality, for any stopping time  $s>0$,  
\beaa
V(x)&=&\sup_{{\pi}(\cd)\in\sA (x)}\hE_x\Big[\int_0^{s\wedge\tau^\pi}e^{-ct} r_\l^\pi(t) dt +e^{-c(s\wedge\tau^\pi)}V(X_{s\wedge\tau^\pi}^\pi)\Big], \q s>0.
\eeaa
%Here in the above $\{\pi_t\}$ is the density process of an open-loop control $\pi(\cd)\in\sA(x)$. 
Noting  that $V(0)=0$, we can (formally) argue  that $V$ satisfies the HJB equation:
\bea
\label{hjbfc}
\qq \q \begin{cases}
\dis  cv(x)\neg=\neg\sup_{{\pi}\in\hL^1[0,a]}\int_0^a\neg\Big[w\neg-\neg\lambda\ln\pi(w)+\frac{1}{2}\sigma^2v''(x)\neg+\neg(\mu-w)v'(x)\Big]\pi(w)dw;\\
v(0)=0.
\end{cases}
\eea
%Here in the above $\hL^1_+([0,a])$ is the set of density function of measures of $\sP([0,a])$.  
Solving the calculus of variation   problem $\sup_{{\pi}\in\hL^1[0, a]}\int_0^a F(x, w,\pi) dw$, where $F(x, w,\pi):= [(1-v'(x))w - \lambda\ln \pi ]\pi $, $\pi\in [0, a]$, and the corresponding Euler-Lagrange equation  $ \frac{\partial}{\partial  \pi}  F (x,w, \pi)= 0$ (see \cite{bff}), we see that the optimal feedback control has the following {\it Gibbs form}:
%. Solving the equation $ \frac{\partial}{\partial y}  F_x (w,y)= 0$ For a fixed $x \in \hR_+$, denote  where  \rightarrow{ \hR}$. Then the maximization problem on the right-hand side can be rewritten as, and multiplying by a normalizing constant using the fact that $\pi \in \sP([0, a])$, we readily obtain that the optimal feedback control has the following {\it Gibbs form}, assuming all derivatives exist:}
% {\color{red}
% \bea
% \label{zykz}
% {\pmb\pi}^*(w,x)=\cfrac{1-v'(x)}{\lambda\left[e^{\frac{a}{\lambda}(1-v'(x))}-1\right]}\cdot e^{\frac{w}{\lambda}(1-v'(x))}\1_{\{v'(x)\neq1\}}+
% \frac{1}{a}\1_{\{ v'(x)=1\}}.
% \eea
% }
\bea
\label{zykz}
\left\{\ba{lll}
{\pmb\pi}^*(w,x)= G(w, 1-v'(x)), \ms\\
\dis G(w,y) := \frac{y}{\lambda [e^{\frac{a}{\lambda}y}-1]}\cdot e^{\frac{w}{\lambda}y}\1_{\{y \neq 0\}}+
\frac{1}{a}\1_{\{ y=0\}}, \q y \in \hR. 
\ea\right.
\eea
%where $$. 
Plugging \eqref{zykz} into (\ref{hjbfc}), we see that the HJB equation (\ref{hjbfc}) becomes a second order ODE:
\bea
\label{ode}
\frac{1}{2}\sigma^2v''(x)+f(v'(x))-cv(x)=0, \q x\ge 0;\qq
v(0)=0,
\eea
where the function $f$ is defined by 
\bea
\label{fdxx}
f(z):=\Big\{ \mu z +\lambda \ln\Big[\frac{\lambda(e^{\frac{a}{\lambda}(1-z)}-1 )}{1-z}\Big]\Big\}\1_{\{z\neq1\}}+[
\mu+\lambda\ln a]\1_{\{z=1\}}.
\eea

The following result regarding the function $f$  is important in our discussion. Since the proof is rather elementary, we omit it and refer the interested reader to \cite{BGMX} for details. 
\begin{prop}
\label{fdxz}
The  function $f$ defined by \eqref{fdxx} enjoys the following  properties:

(1) $f(z)= \mu z + \l \ln w(z) $ for all $z \in \hR$,  where $w(z) = a+ \sum_{n=2}^\infty \frac{ a^n (1-z)^{n-1}}{n! \l^{n-1}}$, $z\in\mathbb{R}$. {\color{black}In particular,   $f \in \hC^\infty(\mathbb{R})$ and is convex on $\hR$;
 
 \ms
 %all $z \in \mathbb{R}$$  is second-order continuously differentiable and convex on $\mathbb R$;\

(2) There exists a unique point $z_0\in  (1,1+\lambda)$ such that
%. the function $f$ is convex and has a unique intersection point with 
$f(z_0)-\mu z_0=\l\ln w(z_0)=0$.
%, $x\in\hR$. Moreover, the abscissa value of the intersection point 
\qed}
\end{prop}

We should note that (\ref{ode}) can be viewed as either a boundary value problem of an elliptic PDE with unbounded domain $[0, \infty)$  or a second-order ODE defined on $[0,\infty)$. But in either case,  there is missing information on boundary/initial conditions. 
Therefore, the well-posedness of the classical solution is actually not trivial.
%we now briefly describe. 
 
Let us first consider the equation (\ref{ode}) as an ODE defined on $[0,\infty)$. Since the value function is non-decreasing by Proposition \ref{assl}, for the sake of argument let us first  consider (\ref{ode}) as an ODE with initial condition
$v(0)=0$ and $v'(0)= \tilde \a>0$. 
By denoting $X_1(x) = v(x) $ and $X_2(x) = v'(x)$, we see that (\ref{ode}) is equivalent to the following system of first-order ODEs: for $x\in[0,\infty)$,
\begin{eqnarray} \label{OODE}
\begin{cases}
X_1' = X_2, \qq\qq\qq\qq &X_1(0)=v(0)=0;\\
X_2'= \frac{2c}{\sigma^2} X_1-\frac{2}{\sigma^2}f(X_2), &X_2(0)=v'(0).
\end{cases}
\end{eqnarray}
Here $f$ is an entire function. Let us define $\tilde X_1:= X_1 - \frac{f(0)}{c}$, $X:=(\tilde X_1, X_2)^T$,  $A := \begin{bmatrix} 0 & 1 \\ \frac{2c}{\sigma^2} & -\frac{2}{\sigma^2} h'(0) \end{bmatrix}$ and $q(X) = \begin{bmatrix} 0 \\ -\frac{2}{\sigma^2} k(X_2) \end{bmatrix}$ where  $h(y): = f(y)-f(0) = yh'(0)+ \sum_{n=2}^\infty \frac{h^{(n)}(0)y^n}{n!} = yh'(0)+k(y) $.
 Then,  $X$ satisfies the following system of ODEs:
\bea
\label{ODEN}
X'=AX+q(X), \qq X(0)=(- f(0)/c, v'(0))^T.%\begin{cases}
%
%\tilde X_1' = X_2;\\
%X_2'= \frac{2c}{\sigma^2} \tilde X_1-\frac{2}{\sigma^2}[h'(0) X_2 +k(X_2)],
%\end{cases}
\eea
%This is of the format $$ where $ X= \begin{bmatrix} \tilde X_1 \\ X_2 \end{bmatrix}, 
It is easy to check $A$ has  eigenvalues  $\l_{1,2} = \frac{-h'(0) \mp \sqrt{2c\sigma^2 + h'(0)^2}}{\sigma^2}$, with
$\l_1<0<\l_2$.
Now, let $Y=PX$, where $P$ is 
such that $PAP^{-1} = \mbox{\rm diag}[ \l_1,  \l_2 ]:= B$. Then $Y$ satisfies 
\bea
 \label{ODEP}
Y' = BY+g(Y),  \q Y(0)=PX(0),
\eea
 where $g(Y) = P q(P^{-1}Y)$. 
Since $\nabla_Yg(Y)$ exists and tends to 0 as $|Y| \to 0$, and  $\l_1<0<\l_2$,  
we can follow \cite[Theorem 13.4.3]{CoLe} to construct a solution 
$\tilde \phi $ to \eqref{ODEP}, 
%  for certain values of $\tilde \a$, 
such that  $|\tilde \phi(x)| \leq C_1e^{-C_2 x} $ for some constants $C_1,C_2>0$. 
Hence $|\tilde \phi(x)|\to 0$, as $ x \to \infty$. Thus, the function  $ \phi(x):=P^{-1} \tilde \phi(x) $ is a solution to \eqref{ODEN}
 satisfying $|\phi(x)|\to0$,  as $x \to \infty$. In other words,   \eqref{OODE} has a solution $(X_1(x), X_2(x)) \to 
 %(0+\frac{f(0)}{c}, 0) = 
 (\frac{f(0)}{c},0) $ as $x \to \infty$. 
 We summarize the discussion above as the following result.  
 
 \begin{prop}
 \label{classic}
The differential equation (\ref{ode}) has a  classical solution  $v$  
 %defined on $[0, \infty)$, satisfying $v'(0)=\a$. Furthermore, $v$ 
that enjoys the following properties:
 
 (i) $v'(0)>0$ and $\lim_{x\to\infty}  ( v(x), v'(x))=(\frac{f(0)}{c},0)$;
 
 \ss
 (ii) $v$ is increasing and concave. 
 % \ms
% (iii)
 \end{prop}
 
 {\it Proof} See Appendix. 
 \qed

%\ss
%{\bf Viscosity Solution of \eqref{ode}.}
We note that Proposition \ref{classic} requires that 
%cannot be applied directly since we do not know 
$v'(0)$ exists, which is not a priori known. We now justify this point briefly using the theory of viscosity solutions. We first consider  (\ref{ode}) as an elliptic PDE defined on $[0,\infty)$, and recall (cf. e.g.\cite{BM}) that a viscosity solution of (\ref{ode}) is said to be of class $(L)$ if it is bounded and increasing on $\hR_+$. Next, let us define
\bea
\label{supsolution}
\underline{\psi}(x):=1-e^{-x},\q 
\overline{\psi}(x):= \frac{A}{M}(1-e^{-M(x\wedge b)}),   \qquad x\in[0,\infty), 
%\frac{A}{M}(1-e^{-Mb}),\qquad x\in[b,\infty),
%\end{cases}
\eea
where $A$, $M$, $b>0$ are constants  satisfying $M>2\mu/\sigma^2$ and the following constraints:
\bea
\label{AMbound}
 \left\{\ba{lll}
\frac{1}{M}\Big\{\ln\Big(\frac{A}{A-M}\Big)\vee \ln\Big(\frac{A}{A-\frac{f(0)}{c}M}\Big)\Big\}< b<\frac{1}{M}\Big\{\ln\frac{A}{H}\wedge \ln\Big(\frac{\sigma^2}{2\mu}M\Big)\Big\}; \ms\\
%\end{eqnarray*}
%and
%\begin{eqnarray*}
 A>\max\left\{M+H,~\frac{f(0)}{c}M+H,~\frac{\sigma^2M^2}{\sigma^2M-2\mu},~\frac{f(0)}{c}\cdot\frac{\sigma^2M^2}{\sigma^2M-2\mu}\right\}. 
\ea\right.
\eea
%\vspace{-2mm}
Then, following the argument of \cite{BGMX} or  \cite{BM}, one can show that
$\underline{\psi}(\cd)$ (resp. $\overline{\psi}(\cd)$) is a viscosity  sub- (resp. super-) solution of \eqref{ode} on $[0,\infty)$, and that  $\underline{\psi}(x)\leq\overline{\psi}(x)$ on $[0,\infty)$.
Now consider the set 
\bea
\label{FF}
\qq\q  \mathfrak{F}:=\{u\in \mathbb{C}(\hR_+)\ | \ \underline{\psi}\le u\le\overline{\psi};~ \mbox{$u$ is a class ($L$) vis. subsolution to \eqref{ode}}\}.
\eea
%Clearly, $\underline{\psi}\in \mathfrak{F}$, so $\mathfrak{F}\neq\emptyset$. Define
%\bea
%\label{udef}
Define $ \hat v(x)=\sup_{u\in\mathfrak{F}}u(x)$, $x\in[0,+\infty)$, and
%\eea
%and let $v^*$ (resp. $v_*$) be the USC (resp. LSC) envelope of $\hat v$, defined respectively by
\beaa
v^*(x) (\mbox{resp. $v_*(x)$}):=\limsup_{r\downarrow0} (\mbox{resp. $\dis\liminf_{r\downarrow0}$})\{\hat v(y):y\in(0,+\infty), |y-x|\leq r\},
\q x\in[0,+\infty).
%&&=\{v(y):y\in(0,+\infty), |y-x|\leq r\}.
\eeaa
Then, similar to \cite{BM}, one shows that $v^*$ (resp. $v_*$) is a viscosity sub-(resp. super-)solution of class $(L)$ to \eqref{ode} on $\hR_+$, hence $v^*\leq v_*$, thanks to Comparison Principle. Since by definition  $v_*\le v^* $,  we conclude that $v^*=v_*=\hat v$ is the unique viscosity solution of class $(L)$ to \eqref{ode} on $\hR_+$, and we shall now simply denote $v=\hat v$.
Finally, note that for small $x>0$ we have
$$ \frac{1-e^{-x}}{x}=\frac{\ul{\psi}(x)}{x}\le \frac{v(x)}{x}\le \frac{\ol{\psi}(x)}{x}=\frac{A}{M} \frac{(1-e^{-Mx})}{x}.
$$
Sending $x\searrow 0$ we obtain that $1\le \liminf_{x\searrow0}\frac{v(x)}{x}\le \limsup_{x\searrow0} \frac{v(x)}{x}\le A$. 
Now, recall from Proposition \ref{classic}   that the ODE (\ref{ode}) has a bounded $\hC^2 $-solution
satisfying $v'(0+)=\tilde \a$, and it is increasing and concave, whence a viscosity solution to (\ref{ode}) of class ($L$). Thus by 
uniqueness we conclude that
%Theorem \ref{Comp}, the bounded viscosity solution to (\ref{ode}) of class ($L$) is unique, thus 
the viscosity solution $v\in\hC^2_b(\hR_+)$ and  $v'(0+)=\tilde \a = \liminf_{x\searrow0}\frac{v(x)}{x}=\limsup_{x\searrow0} \frac{v(x)}{x}$. 

Having argued the well-posedness of ODE (\ref{ode}) from both the classical and viscosity senses, we now look at its connection to the value function and the optimal strategy. In fact, it is more or less standard to argue (see, e.g., \cite{YongZhou}) that the value function $V$ is a viscosity solution satisfying 
\bea
\label{Vsupv}
V(x)=\sup_{u\in\mathfrak{F}}u(x)=  v(x), \qq x\in[0,+\infty),
\eea
where $\mathfrak{F}$ is defined by (\ref{FF}). 
It then follows from the previous argument that $V'(0+)$ exists and $V$  is the (unique) classical solution of class ($L$) to \eqref{ode}.
Next, consider the feedback control:
\bea
\label{zycl}
\pi^*_t(w)= G(w,1-V'(X^{\pi^*}_t)), 
%\end{cases}
\eea
where $G(\cd,\cd)$ is the Gibbs function defined by (\ref{zykz}). Note that $|r_\l^{\pi^*}(t)|
=\left|\int_{0}^{a}\bar{f}(V'(X_t^*))\pi_t^*(w)dw\right|$, where $f(\cd)$ is defined by (\ref{fdxx}), and it is easy to check that
%$\bar{f}(z):=
%\big\{wz+\lambda \ln\big[\frac{\lambda (e^{\frac{a}{\lambda}(1-z)}-1 )}{1-z}\big]\big\}\1_{\{z\neq1\}}+[w+\lambda\ln a]\1_{\{z=1\}}$. 
%Thus
\beaa
\hE_x\Big[\int_{0}^{\tau^{\pi^*}} e^{-ct}|r_\l^{\pi^*}(t)|dt\Big] 
=\hE_x\Big[\int_0^{\t^{\pi^*}}\neg\neg\neg e^{-ct}\Big|\int_{0}^{a}\bar{f}(V'(X_t^*))\pi_t^*(w)dw\Big|dt\Big]<+\infty,
\eeaa
as $V'(X_t^*)\in(0,V'(0+)]$, thanks to the concavity of $V$. Consequently $\pi^*\in\mathcal{A}(x)$. Since 
%Finally, since $V\in\hC^2_b(\hR_+)$ and 
$\pi^*$ is defined by (\ref{zycl}) and is obviously the maximizer of the Hamiltonian in the HJB equation
(\ref{hjbfc}), a simple application of It\^o's formula then leads to the optimality of $\pi^*$.

\section{Policy Update}
\label{POU}
\setcounter{equation}{0}

We now turn to an important step in the RL scheme, that is, the so-called {\it Policy Update}. More precisely, we prove a {\it Policy Improvement Theorem}  which states that for any closed-loop policy $\pmb\pi\in \cA_{cl}(x)$, we can construct another  $\tilde{\pmb\pi}\in \cA_{cl}(x)$, such that $J(x, \tilde{\pmb\pi})\ge J(x, \pmb\pi)$. Furthermore, we argue that such a policy 
 policy updating procedure can be constructed without using the system parameters, and we shall discuss the convergence of the iterations to the optimal policy.

To begin with, for $x\in\hR$ and
%We consider the open-loop control strategy $\pi$ induced from
$\pmb\pi \in \sA_{cl}(x)$, let $X^{\pmb\pi, x}$ be the unique strong solution to the SDE (\ref{sde2}). For  $t>0$, we consider
the process $\hat W_s:=W_{s+t}-W_t$, $s>0$. Then $\hat W$ is an $\hat\hF$-Brownian motion, where $\hat\cF_s=\cF_{s+t}$, $s>0$. Since the SDE (\ref{sde2}) is time-homogeneous, the path-wise uniqueness then renders the {\it flow property}: $X_{r+t}^{\pmb\pi, x}=\hat X_r^{\pmb\pi,X_t^{\pmb\pi,x} }$, $r\ge 0$, where $\hat X$ satisfies the SDE
\bea
\label{SDE0}
d\hat X_s=\Big(\mu-\int_{0}^{a}w \pmb{\pi} (w, \hat X_s)dw\Big)ds+\sigma d\hat W_s, \q s\ge 0; \q \hat X_0= X_t^{\pmb\pi,x}.
\eea

Now we denote $\hat \pi:=\pmb\pi(\cd, \hat X_\cd)\in\sA_{ol}(X_t^{\pmb\pi, x})$ to be the open-loop  strategy  induced by the closed-loop control $\pmb\pi$.
Then the corresponding cost functional can be written as (denoting $X^{\pmb\pi}=X^{\pmb\pi,x}$)
%given the initial condition $X_0^{\hat \pi} = X_t^{\pmb\pi,x}$; and denote
%with respect to the policy $\pmb\pi \in\cA_{cl}$ and initial condition $X_0= X_t^{\pi} $  by
%$J (X_t^{\pi} ,\pmb\pi)$. Namely
\bea
\label{Jpmbpi}
 \q J (X_t^{\pmb\pi};\pmb\pi) = \hE _{X_t^{\pmb\pi}} \Big[ \int_{0}^{\tau ^ {\pmb{\pi}}_{X_t^{\pmb\pi}}}  e^{-cr}
\Big[\int_{0}^{a}
( w - \lambda \ln \hat{\pi}_r(w)) \hat{\pi}_r(w) dw \Big] dr\Big], \q t\ge 0,
\eea
 where
$ {\tau ^ {\pmb\pi}_{X_t^{\pmb\pi, x}}} = \inf \{r>0: \hat X_r^{\pmb\pi,X_t^{\pmb\pi,x} } <0 \}$.
It is clear that, by the flow property, we have  $ {\tau ^{\pmb\pi}_x} = {\t^{\pmb\pi}_{X_t^{\pmb\pi,x}}}+t$, $\hP$-a.s. on $\{\t^{\pmb\pi}_x>t\}$.
Next, for any admissible policy  $ \pmb{\pi}\in \sA_{cl}$, we formally define a new feedback control policy as follows:
for $(w,x)\in [0,a]\times  \hR^+$, 
\bea
\label{pitilde}
 \pmb{\tilde \pi}(w,x)  :=  G(w, 1- J^{'}(x;\pmb{\pi}) ),
   \eea
where $G(\cd, \cd)$ is the Gibbs function defined by (\ref{zykz}). 
We would like to emphasize that the new policy $\pmb{\tilde\pi}$ in (\ref{pitilde}) 
depends  on $J$ and $\pmb{\pi}$, but is
independent of the coefficients $(\m, \si)$.  To facilitate the argument, we introduce the following definition.
\begin{defn}
\label{strongadm}
A function $x\mapsto {\pmb \pi}(\cd;x)\in\sP([0,a])$  is called ``Strongly Admissible" if  its density function  enjoys the following properties:

\ss
(i) there exist $u,l > 0$ such that $l \leq \pmb \pi (x,w) \leq u$, $x \in \mathbb{R}^+ $ and $ w \in [0,a] $;

\ss
(ii)  there exists $K > 0$ such that $|\pmb \pi(x,w) - \pmb \pi (y,w)| \leq K |x-y|$, $x,y \in \mathbb{R}^+ $, uniformly in $w$.

The set of strongly admissible controls is denoted by $\sA^s_{cl}$. 
\qed
\end{defn}

{\color{black}We are particularly interested in the following type of feedback controls ${\pmb \pi}(\cd;x)\in\sP([0,a])$  whose density takes the form $\pmb\pi(w,x)  =  G(w,c(x))$ where $c\in \hC^1_b(\hR_+)$, and $G$ is the Gibbs function defined by (\ref{zykz}). Since 
%let $K>0$ be such that  $|c(x)|+|c'(x)|\le K$, $x\in\hR_+$. Next, note that 
$G$ is  positive, continuous, and $ G(w,\cd) \in\hC^\infty(\hR)$ for fixed $w$, it is easy to check that  
%we have $  l \le \pmb\pi(x,w)= G(w,c(x))    \le u$, $ (w,x) \in  [0,a] \times \hR^+ $, and  is uniformly Lipschitz on $\hR^+$, uniformly in $w\in[0,a]$,
$\pmb\pi(\cd,\cd) = G(\cd,c(\cd))\in \sA^s_{cl} $ for any $c\in \hC^1_b(\hR_+)$. We state this fact as the following lemma for ready reference.} 
%The following lemma justifies the Definition \ref{strongadm}.
% and Assumption \ref{assumAP} is the following lemma.
\begin{lem}
\label{lemma4}
Suppose that a function $x\mapsto {\pmb \pi}(\cd;x)\in\sP([0,a])$ whose density takes the  form $\pmb\pi(w,x)  =  G(w,c(x))$ where $c\in \hC^1_b(\hR_+)$.
% is continuously differentiable and $ \exists u_1, u_2 > 0$ such that $ |c(x)| \leq u_1 $ and $|c'(x)|\leq u_2 $, 
Then $\pmb\pi\in\sA^s_{cl}$. 
\end{lem}
%{\it Proof.} 
%
%\qed

In what follows, {\color{black} we fix $\l>0$. 
%shall use the following notations. 
For $\pmb\pi\in\sA_{cl}$, we denote}
\bea
\label{ABr}
\q r^{\pmb\pi}(x) =r^{\pmb\pi}_\l(x):= \int_0^a (w -\l  \ln \pmb\pi(w,x)) \pmb\pi(w,x)  dw; \q b^{\pmb\pi}(x)=\m- \int_0^a w \pmb\pi(w,x)  dw.
\eea
Clearly, for $\pmb\pi\in\sA_{cl}^s$,   $ b^{\pmb \pi} $ and $r^{\pmb \pi} $ are bounded and are Lipschitz continuous.  We denote
$X:=X^{\pmb\pi, x}$ to be the solution to SDE (\ref{sde2}), and rewrite the cost function (\ref{Jxpi}) as
\bea
\label{Jxpi1}
\qq J(x, \pmb\pi)=\hE_x\Big[\int_0^{\t_x^{\pmb\pi}} e^{-cs}r^{\pmb\pi}(X^{\pmb\pi, x}_s)ds\Big].
\eea
where $ \t^{\pmb\pi}_x = \inf \{t>0:  X^{\pmb\pi,x}_t <0 \}$. 
 Thus, in light of the Feynman-Kac formula,  for any $\pmb\pi \in \sA^s_{cl}$, $J(\cd, \pmb\pi)$ is the {\it probabilistic solution} to the following ODE on $\hR_+$:
  \bea
\label{fkpde}
\qq L^{\pmb\pi}[u](x)+r^{\pmb\pi}(x)\neg:= \neg\frac{1}{2} \sigma^2 u_{xx}(x)  + b^{\pmb\pi}(x)u_x(x) -cu(x) + r^{\pmb\pi}(x) = 0,   
~ u(0) =0.
\eea 
%Now recall from Proposition \ref{assl} that the value function $V$ is bounded. It follows that $J(\cd, \pmb\pi)\le V(\cd)$ is also uniformly bounded, uniform in $\pmb\pi$. Furthermore, if we 
Now let us denote $u^{\pmb\pi}_R$ to be solution  to the linear elliptic equation
%Finally, we note that as a bounded solution to 
(\ref{fkpde}) on finite interval $[0,R]$ with boundary conditions $u(0)=0$ and $u(R)=J(R, \pmb\pi)$, then by the regularity and the boundedness of $b^{\pmb\pi} $ and $r^{\pmb\pi} $,
%as the solution to (\ref{fkpde}) we should have $J(\cd, \pmb{\pi})\in \hC^2(\hR^+)$. 
and using only the interior type Schauder estimates (cf. \cite{GT}), one can show that $u^{\pmb\pi}_R\in\hC^2_b([0,R])$  and the bounds of $(u^{\pmb\pi}_R)'$ and $(u^{\pmb\pi}_R)''$ depend only on those of the coefficients $b^{\pmb\pi}$, $r^{\pmb\pi}$ and $J(\cd, \pmb\pi)$, but uniform in $R>0$. By sending $R\to\infty$ and applying the standard diagonalization argument (cf. e.g., \cite{Lady}) one shows that $\lim_{R\to\infty}u^{\pmb\pi}_R(\cd)=J(\cd, \pmb\pi)$, which satisfies (\ref{fkpde}). We summarize the above discussion as the following proposition for ready reference. 
%Therefore, for the sake of argument in what follows we shall make use of the following assumption.
\begin{prop}
\label{propDB}
If $ \pmb{\pi} \in\sA^s_{cl}$, then $J(\cd, \pmb{\pi})\in \hC^2_b(\hR^+)$, and the bounds 
of $J'$  and $J''$ depend  only on those of  $b^{\pmb\pi}$, $r^{\pmb\pi}$, and  $J(\cd, \pmb\pi)$. 
\qed
\end{prop}

Our main result of this section is the following {\it Policy Improvement Theorem}.

\begin{thm}
\label{PIthm}
Assume that Assumption \ref{jc}  is in force. Then, let $\pmb\pi\in\sA^s_{cl}$ and let $\tilde{\pmb\pi}$ be defined by (\ref{pitilde}) associate to $\pmb\pi$, it holds that
%Then it holds that
$ J(x,\pmb{ \tilde\pi}  ) \geq J(x,\pmb{ \pi})$,  $x \in  \hR_+$.
\end{thm}

{\it Proof.} Let $\pmb\pi\in\sA^s_{cl}$ be given, and let $ \pmb{\tilde\pi  } $ be the corresponding control defined by (\ref{pitilde}).  
Since $\pmb\pi\in\sA^s_{cl}$, $ b^{\pmb\pi}$ and $r^{\pmb\pi}$ are uniformly bounded, and by Proposition \ref{propDB}, 
%pis finite and as a consequence so is $  \|J(\cd, \pmb\pi) \|_{L^\infty(\hR_+)}$. Hence, 
$(1-J'(\cd, \pmb \pi)) \in \hC^1_b(\hR_+)$. Thus  Lemma \ref{lemma4} (with $c(x) = 1 - J'(x , \pmb\pi)$) implies that
$\tilde{\pmb\pi}\in\sA^s_{cl}$ as well. Moreover, since $\pmb\pi\in\sA^s_{cl}$, 
%Moreover, note that for the given $\pmb\pi\in\sA^s_{cl}$, the cost functional 
$J(\cd, \pmb\pi)$ is a $\hC^2$-solution to the ODE (\ref{fkpde}).  
Now recall that $\tilde{\pmb\pi}\in\sA^s_{cl}$ is the maximizer of $\sup_{\h{\pmb\pi}\in\sA^s_{cl}} [b^{\h{\pmb\pi}}(x)J'(x, \pmb\pi)+r^{\h{\pmb\pi}}(x)]$,  we have
\bea
\label{Jpitilde}
L^{\tilde{\pmb\pi}}[J(\cd, \pmb\pi)](x)+ r^{\tilde{\pmb\pi}}(x) \ge 0, \q x\in\hR_+.
\eea
Now, let us consider the process $X^{\tilde{\pmb\pi}}$, the solution to  (\ref{sde2}) with $\pmb\pi$ being replaced by $\tilde{\pmb{\pi}}$. Applying It\^o's formula
to  $ e^{-ct}J(X^{\tilde{\pmb\pi}}_t, {\pmb{\pi}}) $  from  $0$ to $\t^{\tilde{\pmb\pi}}_x\wedge T$, for any $T>0$, and noting the definitions of $b^{\tilde{\pmb\pi}}$ and $r^{\tilde{\pmb\pi}}$, we deduce from (\ref{Jpitilde}) that
%under the state process, we have
%for $  t \in [0, \tau^{\tilde\pi}_x]$
$$ e^{-c(\t^{\tilde{\pmb\pi}}_x\wedge T)}J (X_{\t^{\tilde{\pmb\pi}}_x\wedge T}^{\tilde{\pmb\pi}} ,\pmb{ \pi}) \ge J (x,\pmb{ \pi})   - \int_0^{\t^{\tilde{\pmb\pi}}_x\wedge T} \neg\neg e^{-cr}r^{\tilde{\pmb\pi}} (X_r^{{\tilde\pi}}) dr +   \int_0^{\t^{\tilde{\pmb\pi}}_x\wedge T} \neg
\neg e^{-cr}  J^{'}(X_r^{\tilde{\pmb\pi}},\pmb{\pi})  {\sigma}  dW_r.
$$
%where the last inequality above is due to (\ref{Jpitilde}). 
Taking expectation on both sides above, sending $T\to \infty$  and noting that $J (X_{\t^{\tilde{\pmb\pi}}_x}^{\tilde{\pmb\pi}} ,\pmb{ \pi})=J(0, \pmb\pi)=0$, we obtain that $J(x,\pmb{ \pi}) \le J(x, \tilde{\pmb\pi})$, 
$x \in \hR^+$, proving the theorem.
\qed

In light of Theorem \ref{PIthm}, we can naturally define a ``learning sequence" as follows.  We start with $c_0\in C^1_b(\hR^+)$, and define  $\pmb\pi_0(x,w):=G(w, c_0(x)) $,  and $v_0(x) := J(x, \pmb\pi_0)$, 
\bea
\label{pin}
\q \q \pmb\pi_n(x,w)  := G(w, 1- J'(x, \pmb\pi_{n-1})), \qq (w,x)\in [0,a]\times  \hR_+ , \text{ for }\q n \geq 1. 
  \eea 
Also for each $n \geq 1$,  let $v_{n}(x) := J(x, \pmb\pi_{n}) $. The natural question is whether this learning sequence is actually a ``maximizing sequence", that is, $v_n(x)\nearrow v(x)$, as $n\to\infty$. Such a result would obviously justify the {\it policy improvement} scheme, and was proved in the LQ case in \cite{zhou2}. 

Before we proceed, we note that by Proposition \ref{propDB}
the learning sequence $v_n =J(\cd, \pmb{\pi}^n)\in \hC^2_b(\hR_+)$, $n \geq 1$, but the bounds may depend on the coefficients
$b^{\pmb\pi^n}$, $r^{\pmb\pi^n}$, thus may not be uniform in $n$. But by definition $b^{\pmb\pi^n}$ and  Proposition \ref{assl}, we see that $\sup_n\|b^{\pmb\pi^n}\|_{\hL^\infty(\hR_+)}+\|V\|_{\hL^\infty(\hR_+)}\le C$ for some $C>0$. Moreover, since for each $n \geq 1$, $ J(\cd, \pmb\pi^0) \le J(\cd, \pmb\pi^n) \le V(\cd)$, if we choose $\pmb\pi^0\in\sA^s_{cl}$ be such that $J(x,\pmb\pi^0)\ge 0$ (e.g., $\pmb\pi^0_t\equiv \frac1a$), then we have $  \|J(\cd, \pmb\pi^n) \|_{L^\infty(\hR_+)}\le\|V\|_{\hL^\infty(\hR_+)}\le C $ for all $n \ge 1$. That is, $v_n$'s are uniformly bounded, and uniformly in $n$, provided that $r^{\pmb\pi^n}$'s are. The following result, based on the recent work \cite{HWZcvg}, is thus crucial.
\begin{prop}
\label{assumAP}
The functions $r^{\pmb\pi_n}$, $n\ge 1$ are uniformly bounded,  uniformly in $n$. Consequently, 
the learning sequence $v_n =J(\cd, \pmb{\pi}^n)\in \hC^2_b(\hR_+)$, $n\ge 1$, and the bounds of $v_n$'s, up to their second derivatives, are uniform in $n$. 
\qed
\end{prop}
Our main result of this section is the following.
\begin{thm}
\label{convergence}
Assume that Assumption \ref{jc}  is in force. Then the sequence $\{v_n\}_{n \geq 0}$ is a maximizing sequence. Furthermore, 
%, converges to the Optimal Value function $V$ 
 the sequence $\{\pi_n\}_{n \geq 0}$ converges to the optimal policy $\pi^*$. 
\end{thm}
{\it Proof.} 
We first observe that by Lemma \ref{lemma4}  the sequence $\{\pmb\pi_n\}\subset \sA_{cl}^s$, provided $\pmb\pi_0\in\sA_{cl}^s$. Since
$v_n=J(\cd, \pmb\pi_n)$,  Proposition \ref{assumAP} guarantees that $v_n\in \hC^2_b(\hR_+)$, and the bounds are independent of $n$; Theorem \ref{PIthm} shows that $\{v_n\}$  is monotonically increasing, and thus  $\{v_n\}_{n\ge0}$ must converge, say, to $v^*(\cd)$.
Let us fix any compact set $E\subset \hR_+$. A simple application of 
Arzella-Ascolli Theorem shows that there exist a subsequence $\{n_k\}_{k\ge1}$ such that $\{v'_{n_k}\}_{k \geq 0}$ converge  uniformly on E, say, to $v^{**}(\cd)$.
Since obviously $\lim_{k\to\infty}v_{n_k} =v^{*} $, noting that the derivative operator is a closed operator, it follows that $v^{**}(x)= (v^*)'(x)$, $x \in E$. By the same argument, for any subsequence of $\{v'_n\}$, there exists a sub-subsequence that converges uniformly on $E$ to the same limit $(v^*)'$, and thus the sequence $\{v'_n\}$ itself converges uniformly on $E$ to $(v^*)'$. Since $E$ is arbitrary, this shows that $\{(v_n, v'_n)\}_{n\ge0}$ converges uniformly on compacts to $(v^*, (v^*)')$. 
 Since $\pmb \pi_n $ is a continuous function of $v'_n$, we see that $\{\pmb \pi_n\}_{n \geq 0}$  converges uniformly to $\pmb\pi^*\in\sA_{cl} $ defined by $\pmb\pi^*(x,w)  := G(w, 1- (v^*)'(x))$. 
%\eea}

Finally, applying Lemma \ref{lemma4} we see that $\pmb\pi^* \in\sA^s_{cl}$, and the structure of the $\pmb\pi^*(\cd,\cd) $ guarantees that $v^*$ satisfies the HJB equation (\ref{hjbfc}) on the compact set $E$. Since $E$ is arbitrary, we see that 
%By expanding the result to  using the fact that, 
$v^*$ satisfies the HJB equation (\ref{hjbfc}) (or equivalently (\ref{ode})) on $\hR_+$. 
Now, by using the slightly modified verification argument in Theorem 4.1 in  \cite{HWZcvg}, we conclude that $v^* = V^*$ is the unique solution to the HJB equation (\ref{hjbfc}) and thus $\pi^*$ by definition is the optimal control. 
\qed

\begin{rem}
\label{remark4.4}
{\rm
In theory, since  Theorem \ref{PIthm} is based on finding the maximizer of the Hamiltonian, the learning strategy (\ref{pin}) may depend on the system parameters, making a model-based RL approach more suitable. However, a closer look at the learning parameters $c_n$ and $d_n$ in (\ref{pin}) shows that, in our context, they depend only on $v_n$, but not $(\m, \si)$ directly, allowing us to follow a model-free RL approach using the learning strategy (\ref{pin}) for our numerical analysis in \S7. 
\qed}
\end{rem}

\section{Policy Evaluation --- A Martingale Approach}
\setcounter{equation}{0}

Having proved the policy improvement theorem, we turn our attention to an equally important issue in the learning process, that is, the {\it evaluation} of the cost (value) functional, or the {\it Policy Evaluation}. The main idea of policy evaluation in the reinforcement learning literature usually refers to a process of approximating the cost functional $J(\cd,\pmb\pi)$, for a given feedback control $\pmb\pi$, by approximating $J(\cd,\pmb\pi)$ with a parametric family of functions $J^\theta $, where $\theta \in \Theta \subseteq \mathbb{R}^l$.
Throughout this section, we shall consider a fixed feedback control policy $\pmb\pi\in\sA^s_{cl}$. Thus, for simplicity of notation, we shall drop the superscript $\pmb\pi$ and write $r(x) = r^{\pmb\pi}(x), 
b(x) = b^{\pmb\pi}(x) $ and $ J(x, \pmb\pi)=J(x)$; $\t_x = \t^{\pmb\pi}_x$. 

  We note that for $\pmb\pi\in\sA^s_{cl}$, the functions $r,b \in \hC_b(\hR_+)$ and $J \in \hC_b^2(\hR_+)$. 
Now let $X^x=X^{\pmb\pi, x}$ be the   solution to the SDE (\ref{sde2}), and 
$J(\cd)$ satisfies the ODE (\ref{fkpde}). Then,  applying It\^o's formula we see that
\bea
\label{M}
M^x_t: = e^{-ct}J(X^x_t) + \int_0^t e^{-cs}r(X^x_s) ds, \qq t\ge 0,
\eea
is an $\hF$-martingale. Furthermore, the following result is more or less standard.
\begin{prop}
\label{tildeM}
%For all $x \in \hR_+$, the process $M^x= \{ M^x_t ; t\ge 0\}$ defined by (\ref{M})  is an $\hF$-martingale.  Conversely,
Assume that Assumption \ref{jc} holds, and suppose that $\hat J(\cd)\in\hC_b(\hR_+)$ is such that $\hat J (0) = 0$,  and for all $x \in \hR_+$, the process $\hat M^x:= \{ \hat M^x_t =  e^{-cs}\hat J(X^x_s) + \int_0^te^{-cs} r(X^x_s) ds ;~ t\geq 0\}$ is an $\hF$-martingale. Then $J \equiv \hat J$. 
\end{prop}
{\it Proof.}  First note that $J(0)=\hat{J}(0)=0$, and $X^x_{\t_x}=0$. By (\ref{M}) and definition of
$\hat M^x$ we have
$ \hat M^x_{\t_x}=\int_0^{\t_x} e^{-cs}r(X^x_s) ds =M^x_{\t_x}$. 
%$$ \hat M^x_{\t_x}=\int_0^{\t_x} e^{-cs}r(X^x_s) ds=M^x_{\t^_x}.$$
Now, since $r$, $J$, and $\hat J$ are bounded,  both $\hat M^x $ and $M^x$ are uniformly integrable $\hF$-martingales,  by optional sampling it holds that
$$\hat J(x)=\hat M^x_0=\hE[\hat M^x_{\t_x}|\cF_0]=\hE[M^x_{\t_x}|\cF_0]=M^x_0 = J(x), \qq x\in\hR_+.
$$
The result follows.
 \qed

{\color{black} 
We now consider a family of functions $\{J^\th(x): (x, \th)\in\hR_+\times \Th\}$, where $\Th\subseteq\hR^l$ is a certain index set.
For the sake of argument, we shall assume further that $\Th$ is compact.
Moreover, we shall make the following assumptions.
% for the parameterized family $\{J^\th\}$.
\begin{assum}
\label{ape1}
(i) For a fixed $\theta$, the mapping $x \mapsto J^\theta(x)$ is sufficiently smooth, so that all the derivatives required exist
in the classical sense.

\ss

% (ii) For all $\theta \in \Theta$, $\varphi^\th(X^x_\cd)$ are square-integrable continuous processes, and the mappings $\th\mapsto \|\varphi^\th\|_{\hL^2_{\cF}([0,T])}$ are continuous, where $\varphi^\th=J^\th, (J^\th)'$.

% functions of $\theta$ and their $L^2$-norms are continuous functions of $\theta$.
 
\ss
(ii) There exists a continuous function $ K(\cd)>0$,   such that $ \|J^\theta\|_\infty +\|(J^\theta)'\|_\infty \leq K(\theta) $.
%,  $x\in\hR_+$.
\qed
\end{assum}

We note that the parameterized functions $J^\th$ are utilized to approximate the functions $\{v_n\}_{n \geq 1}$ defined in \S6, and the assumptions in Assumption \ref{ape1} are hence justified by Proposition \ref{assumAP}.

 We now follow the idea of \cite{zhou} to construct the so-called {\it Martingale Loss Function}.
For any $\th\in\Th$, consider the parametrized approximation of the process $M=M^x$:
\bea
\label{Mth}
M^\th_t=M^{\th, x}_t:=e^{-ct}J^\th(X^x_t)+\int_0^t e^{-cs}r(X^x_s)ds, \qq t \geq 0
\eea
 To simplify presentation, we shall
 %In what follows, we shall often 
 drop the superscript $x$ from the processes $X^x$, $M^x$ etc., if there is no danger of confusion.  
In light of Proposition \ref{tildeM}, the best approximation  $J^\th$ is one whose corresponding approximating process $M^\th$ defined by (\ref{Mth}) is a martingale (in which case $J^\th=J$(!)). We define the {\it Martingale Loss function}  
\bea
\label{ML}
\q {ML}(\theta)\neg=\neg \frac{1}{2}  \hE \Big[\neg \int_0^{\t_x}\neg\neg  |M_{\t_x}\neg\neg -\neg {M}_t^\theta|^2 dt\Big]\neg=\neg
 \frac{1}{2}  \hE \Big[\neg \int_0^{\t_x}\neg \neg\neg \big|e^{-ct} J^\theta(X_t)\neg -\neg\neg \int_{t}^{\t_x}\neg\neg \neg e^{-cs} {r} (X_s) ds\big| ^2 dt\Big]. %\nonumber
\eea
We note that the last equality above indicates that the martingale loss function is actually independent of the function $J$, which is one of the main features of this algorithm.
 Furthermore, we define the {\it Discounted  Mean-Squared Value Error (DMSVE)},
\bea
\label{DMSVE}
{\it DMSVE}(\theta) &=& \frac{1}{2}  \hE\Big[ \int_0^{  \t_x} e^{-2ct}  |J^\theta(X_t) - J(X_t) |^2 dt\Big].
 \eea
The following result connects the minimizers of $ML(\cd)$ and $DMSVE(\cd)$. 
% the martingale loss function and the mean-squared error.
{\color{black}
\begin{thm}
\label{equiv}
Assume that Assumption \ref{ape1} is in force. Then, it holds that
\bea
\label{equivLoss}
\arg \min_ {\theta \in \Theta} {\it ML}(\theta) = \arg \min_ {\theta \in \Theta} {\it DMSVE}(\theta).
\eea
\end{thm}
{\it Proof.}  
By definition (\ref{ML})   we can easily rewrite
\bea
\label{2ML}
2{\it ML}(\theta)
%\neg &\neg=\neg&\neg  \hE \Big[ \int_0^{ \t_x} |M_{\t_x}\neg -\neg M_t^\theta|^2 dt\Big]
% = \hE\Big[ \int_0^{\tilde\t_x} |M_{\t_x} - M_t^\theta|^2 dt \Big]  \\
%&=& \hE \Big[ \int_0^{T} |M_{\t_x} -\tilde M_t^\theta|^2 dt \1_{\{{\t_x}\ge T\}}
%+ \int_0^{T} |M_{\t_x} -\tilde M_t^\theta|^2 dt \1_{\{{\t_x} < T\}}\Big]\nonumber \\
%= \hE \Big[ \int_0^{  \t_x} |M_{\t_x}\neg -M_t+M_t  -\neg M_t^\theta|^2 dt\Big]\nonumber \\
%\neg &\neg=\neg&\neg 
=\hE \Big[ \int_0^{   \t_x} \big[ |M_{\t_x}\neg -\neg   M_t|^2\neg +\neg |   M_t \neg- M_t^\theta|^2 +2(M_{\t_x}\neg -  M_t)( M_t - M_t^\theta)  \big] dt\Big]. 
\eea
Noting  (\ref{M}) and (\ref{Mth}), 
%and again the fact $\tilde M^{\th}_t  = \tilde M_t = M_{\t_x}$, for $t \in (\tilde \t_x , T)$, 
we see that
\bea
\label{MMth}
\hE \Big[\neg \int_0^{  \t_x}\neg | M_t -  M_t^\theta|^2  dt\Big] =\hE\Big[\neg \int_0^{   \t_x}\neg\neg e^{-2ct} |J(X_t)  - J^\theta(X_t) |^2  dt\Big] = 2\text{DMSVE}(\theta).
\eea

On the other hand, since $M$ is uniformly integrable, by 
%for a fixed $t \geq 0$, Doob's
optional sampling theorem we have $\hE[M_{\t_x\vee t}|\cF_t]=M_t$, $t\in[0,T]$. Furthermore, note that
 $M_{\t_x\vee t}=M_{\t_x}$ on $\{\t_x\ge t\}\in\cF_t$, we have
\bea
\label{ostep}
\hE\big[(M_{\t_x}-M_t)\1_{\{\t_x\ge t\}}\big|\cF_t\big]
=\1_{\{\t_x\ge t\}}\big(\hE[M_{\t_x\vee t}|\cF_t]-M_t\big)=0,\qq t\ge0 .
\eea
Moreover, since $r, J$ and $J(\theta)$ are bounded, a simple application of 
%(by a coefficient possibly depending on $\theta$, we have that $|M_{\t_x}-M_t|\,|M_t-M^\th_t| \leq e^{-2ct}\Lambda(\theta)$ for all $t \geq 0$ for some $\Lambda(\theta)>0$. Thus, we have
%\beaa
%\int_0^\infty\hE\Big[|M_{\t_x}-M_t|\,|M_t-M^\th_t|\1_{\{\t_x\ge t\}}\Big]dt
%\le \Lambda(\theta) \int_0^\infty e^{-2ct}dt<\infty.
%\eeaa
%Thus, 
Fubini's theorem and (\ref{ostep}), leads to that
\bea
\label{crossterm}
&&  \hE\Big[\int_0^{\t_x}(M_{\t_x}-M_t)(M_t-M^\th_t)dt\Big]= \hE\Big[\int_0^{\infty}(M_{\t_x}-M_t)(M_t-M^\th_t) \1_{\{\t_x\ge t\}}dt\Big]  \\
\neg\neg&\neg\neg=\neg\neg&\neg\neg \int_0^\infty\hE\Big[(M_{\t_x}-M_t)
(M_t-M^\th_t) \1_{\{\t_x\ge t\}}\Big]dt =\neg\neg \int_0^\infty\hE\Big[\hE\big[(M_{\t_x}-M_t)\1_{\{\t_x\ge t\}}\big|\cF_t\big]
(M_t-M^\th_t)\Big]dt=0 .\nonumber
\eea
Combining (\ref{2ML}), (\ref{MMth}) and (\ref{crossterm}), we have 
\bea
\label{MLsplit}
2{\it ML}(\th)=2\,{\it DMSVE}(\th)
+\hE\Big[\int_0^{\t_x}|M_{\t_x}-M_t|^2dt\Big].
\eea
Finally, note that $\hE [ \int_0^{\t_x}  |M_{\t_x} -  M_t|^2   dt]$  is independent of $\theta$,    the result follows.
\qed

%\begin{rem}
%\label{remark5}
We remark that Theorem \ref{equiv} amounts to saying that the best parametrization of $J^\th$  in the sense of DMSVE actually is the minimizer of the martingale loss. Thus, in the rest of the paper, especially in our algorithm design, we shall look for the optimal $\th^*$ in terms of DMSVE. 
%
% $\theta^*$ is a minimizer of  the Martingale Loss functional $ML(\cd)$, it minimizes $DMSVE(\cd)$, and thus, $J^{\theta^*}$ would be an acceptable approximation of $J$. In the rest of the section, we shall therefore focus on the identification of $ \ theta^*$.
%\end{rem}

To end this section, we analyze two necessary approximations that are needed in our algorithm. First, since the ruin time $\t_x$ is unbounded,  our problem is essentially an infinite horizon problem. But on the other hand, our algorithm will be based on the idea of ``batch learning", which requires the simulation of many paths of the state dynamics $X$, thus 
%and involves the random, possibly infinite, interval $[0,\t_x]$, 
any implementable algorithm can only operate on a fixed finite horizon with a finite time grid. We therefore  introduce a {\it time-truncation}
which will restrict our algorithm to a finite time interval $[0,T]$ with a sufficiently large $T>0$. 
%\begin{rem}[Truncation error Analysis]
 %   \label{trunc}

Second, as with any algorithm, we will need to discretize our underlying state process $X$ on a finite time grid. More precisely, let us fix a truncated time horizon $T>0$, and
% To this end, we fix the time horizon $T>0$ (See Remark \ref{trunc} in \S6 for the analysis of the errors of analyzing the truncated problem). We
denote the effective truncated time horizon by $\tilde \t_x:= \t_x\wedge T$. We define the ``truncated" version of the value function and the Martingale Loss function by
 \bea
 \label{vT}
  V_T(x)&:=& \sup_{ \ppi\in\sA^s_{cl} }\hE_x\Big\{\int_{0}^{\tilde \tau_x}e^{-ct}r(X_t) dt\Big\} \qq x\ge 0; \\ 
% \label{MLT}
% {\it ML}_T(\th)&:=&\frac12\,\Big[ \int_0^{\tilde  \t_x} |M_{\t_x}\neg -\neg M_t^\theta|^2 dt\Big].
% % = \frac12\hE\Big[\int_0^{\tilde \t_x}|m(t,\th)|^2dt\Big].
% \eea
%% 
%% define the (truncated) Martingale Loss function,
%\bea
\label{ML}
{\it ML}_T(\th)&:= &\frac12\hE\Big[\int_0^{\tilde \t_x}|e^{-ct}J^\th(X_t)-\int_t^{\tilde \t_x}e^{-cs}r(X_s)\,ds|^2dt\Big].
\eea
% where $\dis m(t,\th):=$.}
%and define the (truncated) Martingale Loss function, consider a partition of $[0, T]$:
{\color{black} We note that in the algorithm, the expectation $\hE^{\hQ}$ over the sampled
actions defined by (2.3)--(2.4) will be implemented as the empirical
average over independently simulated trajectories,
%along each of which the actions are 
sampled from the same fixed policy $\pmb{\pi}$.}

Specifically, let $0=t_0<\cds<t_n=T$ be a partition of
$[0,T]$ with $\D t=t_i-t_{i-1}$, $i=1,\cds,n$.   For $x\in\hR_+$ set 
$ K_x=\min\{  l\in\hN:\exists t \in [l\D t,(l+1)\D t]  :X^x_t<0\}$ so that $\t_x\in[K_x\D t,(K_x+1)\D t]$ and  $\lfloor\t_x\rfloor:=K_x\D t$. We define $N_x=\min\{K_x,n\}$, and hence $t_{N_x}=\lfloor\t_x\rfloor\wedge T$;
both $K_x$ and $N_x$ are integer-valued random variables, and we drop the
subscript $x$ when the context is clear.  In light of (\ref{ML}), the natural  approximation
of the truncated Martingale Loss function $ML_T(\th)$  is then
\bea
\label{MLtildeDt}
2{\it ML}_{T,\D t}(\th)
&:=&\hE\Big[\sum_{i=0}^{N-1}\Big|e^{-ct_i}J^\th(X_{t_i})
-\sum_{j=i}^{N-1}e^{-ct_j}r(X_{t_j})\D t\Big|^2\D t\Big]\\
&=&\hE^{\hQ}\Big[\sum_{i=0}^{N-1}
\Big|e^{-ct_i}J^\th(\tilde X_{t_i})-\sum_{j=i}^{N-1}e^{-ct_j}
\big(a^{\pmb\pi}_{t_j}-\l\ln\pmb\pi(a^{\pmb\pi}_{t_j},\tilde X_{t_j})\big)\D t
\Big|^2\D t\Big], \nonumber
\eea
where in the second equality above, $a^{\pmb\pi}_t$ denotes the action sampled from $\pmb\pi$ at time $t$,   $\tilde X$ the corresponding (simulated) trajectory, and $\hQ=\hQ_T$ is the sampling prior defined by (\ref{Qt}). 
Clearly, the truncated martingale loss is 
%only on the simulated values $\{X_{t_i}\}_{i=0}^{n}$ and is
therefore directly implementable.  
%Under $\hQ$ it takes the form
%\beaa
%2{\it ML}_{T,\D t}(\th)\eeaa
%where 

Now, suppose that for a fixed truncation time $T>0$ and the discretization time step $\D t>0$, we can find an approximation of the minimizer of the truncated martingale loss $ML_{T, \D t}(\cd)$, whose approximation error will conceivably depend on $T$ and $\D t$. We now take a close look at such dependence.  We first give an error bound with respect to $\D t$.
\begin{thm}
\label{conv}
Let Assumptions \ref{jc} and \ref{ape1} be in force and let $T>0$ be fixed, and $\D t \leq 1$.
Then there exists a continuous function $\Gamma(\th)$  such that 
\beaa
\big|{{\it ML}}_{T,\D t}(\th)-{\it ML}_T(\th)\big| \leq \Gamma(\th) \sqrt{\D t}.
\eeaa
That is, $\lim_{\D t \to 0} {{\it ML}}_{T,\D t}(\th) = {\it ML}_T(\th) $ uniformly in $\th$ on compacta.
\end{thm}

{\it Proof.}  See Appendix.
\qed

We next look at the truncation error between $V$ and $V_T$.  We have the following result.
\begin{thm}
\label{Truncerror} 
Under the same assumption as Theorem \ref{conv}, there exists $R>0$, independent of $T$, such that for any $T>0$, it holds that
\bea
\label{VTest}
|V(x)-V_T(x)|\le \frac{R}{c}e^{-cT}, \qq\mbox{ uniformly in $x\ge0$.}
\eea
\end{thm}

{\it Proof.} Recall that the optimal   control  $\ppi^*$ satisfying $V(x)=J(x, \ppi^*) =J^{\ppi^*}(x)$ is given by \eqref{zycl}.  
On the other hand, by Proposition \ref{assumAP},  the approximating sequence $r^{\pmb\pi_n}$, $n\ge 1$, generated by the PIA scheme $\{\pmb\pi_n \}_{n \geq 1}$  are uniformly bounded,  uniformly in $n$, say by a constant $R>0$. Thus we may restrict
the maximization domain in \eqref{vT} to the set
%in both \eqref{v_dup} and  to
$\sA^s_{cl}(R):=\big\{\ppi\in\sA^s_{cl}:\ |r^\ppi(\cdot)|\le R\big\}$,
which contains $\ppi^*$ and the entire sequence $\{\ppi_n\}_{n\ge1}$. 
%Writing
%\bea
%J^\ppi(x):=\hE_x\Big\{\int_{0}^{\t_x^\ppi}e^{-ct}r^\ppi(X_t)\,dt\Big\},
%\qq
%J_T^\ppi(x):=\hE_x\Big\{\int_{0}^{\tilde\t_x^\ppi}e^{-ct}r^\ppi(X_t)\,dt\Big\},
%\eea

Since for every $\ppi\in\sA^s_{cl}(R)$, $\tilde\t_x^\ppi=\t_x^\ppi\wedge T\le\t_x^\ppi$ and  
the two integrands agree on $[0,\tilde\t_x^\ppi]$, we have 
\beaa
\label{perpolicy}
\big|J^\ppi(x)-J_T^\ppi(x)\big|
\;=\;\Big|\hE_x\Big\{\int_{\tilde\t_x^\ppi}^{\t_x^\ppi}e^{-ct}r^\ppi(X_t)\,dt\Big\}\Big|
\;\le\; R\int_{T}^{\infty}e^{-ct}\,dt\;=\;\frac{R}{c}\,e^{-cT}.
\eeaa
Moreover, since $V_T(x)\ge J_T^{\ppi^*}(x)$ and $V(x)=J^{\ppi^*}(x)$,
\bea
\label{VVT}
V(x)-V_T(x)\;\le\; J^{\ppi^*}(x)-J_T^{\ppi^*}(x)\;\le\;\frac{R}{c}\,e^{-cT}.
\eea
On the other hand, we have 
\bea
\label{VTV}
V_T(x)-V(x) =\sup_{\ppi\in\sA^s_{cl}(R)}J_T^\ppi(x)-\sup_{\ppi\in\sA^s_{cl}(R)}J^\ppi(x)
\;\le\;\sup_{\ppi\in\sA^s_{cl}(R)}\big(J_T^\ppi(x)-J^\ppi(x)\big)
\;\le\;\frac{R}{c}\,e^{-cT}.
\eea
Combining (\ref{VVT}) and (\ref{VTV}) we obtain ({\ref{VTest}).
\qed

\begin{rem}
\label{remark5}
By Theorem \ref{Truncerror},  the truncation error (to finite horizon $[0,T]$) is $O(e^{-cT})$, while Theorem \ref{conv} indicates that the time-discretization error is  $O((\D t)^{\frac12})$. If we choose
\bea
\label{Trule}
T \geq \frac{1}{2c}\,\ln\frac{1}{\D t}
\eea
then  $e^{-cT}\le(\D t)^{\frac12}$. That is, the truncation error is dominated by the discretization error. 
%, justifying the truncated time horizon. 
In our actual algorithm design in the next section, we shall choose the truncation time horizon $T$ 
according to \eqref{Trule}, after the discretization size $\D t$ and the
discounting factor $c$ are given.
\qed
\end{rem}
%}

%\bea
%\label{MLT}
%{\it ML}_T(\th):=\frac12\,\Big[ \int_0^{\tilde  \t_x} |M_{\t_x}\neg -\neg M_t^\theta|^2 dt\Big] = \frac12\hE\Big[\int_0^{\tilde \t_x}|m(t,\th)|^2dt\Big].
%\eea
%where $m(t,\th):=e^{-ct}J^\th(X_t)-\int_t^{\tilde \t_x}e^{-cs}r(X_s)\,ds$. Throughout this paper, we use the (truncated) Martingale Loss function defined by \eqref{ML} that depends on the trajectory only through
%$\{X_t\}_{t\le \tilde \t_x }$.  
%
}
}

{\color{black}
\section{Numerical Experiments}
\label{NA}
\setcounter{equation}{0}

In this section, we present a numerical algorithm that reflects the main features of the PIA scheme discussed in the previous sections. 
% Based on the discussions in previous sections, the most suitable method we have available for PE to be used along with Policy Update based on Policy Improvement Theorem is $CTD(\gamma)$ methods. 
%In particular, we shall consider the \textit{CTD}$(0)$ method with special parametrization based on the knowledge of the explicit solution of the original optimal dividend problem (with $\l=0$), but without specifying the market parameter $\m$ and $\si$.
To begin with, we recall the entropy-regularized stochastic control problem (\ref{v}) with (closed-loop) controlled dynamics (\ref{sde2}), and the corresponding HJB equation (\ref{hjbfc}). 
%where
%\bea
%\label{Xpi} 
%dX^\pi_t = \big(\m -\int_0^a w\pi(w, X^\pi_t)dw\big)dt +\si dW_t, \q  t \ge 0; ~ X^\pi_0 = x,
%\eea
%and the cost functional and value function is defined by:
%\bea
%\label{V}
% V (x) = \sup_{\pi} \hE_x\Big[\int_0^{\t_\pi} e^{-ct}\int_0^a (w -\l \ln\pi_t(w))dwdt\Big].
%\eea
%Thus the corresponding is
% \bea
%\label{hjbfc}
%\qq \q \begin{cases}
%\dis  cv(x)\neg=\frac{1}{2}\sigma^2v''(x)\neg+\neg\sup_{{\pi}\in\hL^1[0,a]}\int_0^a\neg\Big[w\neg-\neg\lambda\ln\pi(w)+\neg(\mu-w)v'(x)\Big]\pi(w)dw;\\
%v(0)=0.
%\end{cases}
%\eea
%\section{Possible state dependence of $\m$ ($\m=\m(x)$)}

In order to show the effectiveness and the potential applicability of the RL approach, in what follows we shall design our algorithm allowing the drift coefficient $\m$, which usually includes the interest and premium rate, to be ``state dependent". That is, we assume that $\m=\m(x)$ so that the standard approach for constant parameter estimation does not apply, and we refer to, e.g., \cite{BP, CGY} for such models in optimal dividend problems. It turns out that most of our analysis before can be applied to such an extension without substantial difficulties, which we now briefly describe.

Note that if $\m=\m(x)$, then 
%, which would make the RL approach more significantly meaningful. 
 the HJB equation (\ref{hjbfc}) will take the form:
\bea
\label{hjbfc1}
\qq \q \begin{cases}
%\dis  \neg=\frac{1}{2}\sigma^2v''(x)\neg+\neg\sup_{{\pi}\in\hL^1[0,a]}\int_0^a\neg\Big[w\neg-\neg\lambda\ln\pi(w)+(\mu(x)-w)v'(x)\Big]\pi(w)dw
\dis cv(x)= \frac{1}{2}\sigma^2v''(x)\neg+\mu(x) v'(x)+\neg\sup_{{\pi}\in\hL^1[0,a]}\int_0^a\neg\Big[(1-v'(x))w\neg-\neg\lambda\ln\pi(w) \Big]\pi(w)dw;
\ss\\
v(0)=0.
\end{cases}
\eea
By solving  the calculus of variation problem $\sup_{{\pi}\in\hL^1[0,a]}\int_0^a F(x, w,\pi) dw$, where $F(x, w,\pi) := [(1-v'(x))w - \lambda\ln \pi ]\pi $, $\pi\in [0,a]$, we obtain the same 
%and the corresponding Euler-Lagrange equation  $ \frac{\partial}{\partial  \pi}  F (x,w, \pi)= 0$ to obtain the optimal feedback control has the same {\it 
Gibbs form: ${\pmb\pi}^*(w,x)= G(w, 1-v'(x))$, where $G(\cd, \cd)$ is defined by (\ref{zykz}), and 
%\bea
%\label{Gibbs}
% G(w,y) := \frac{y}{\lambda [e^{\frac{a}{\lambda}y}-1]}\cdot e^{\frac{w}{\lambda}y}\1_{\{y \neq 0\}}+
%\frac{1}{a}\1_{\{ y=0\}}, \q y \in \hR. 
%\eea
%%where $$. 
 the HJB equation (\ref{hjbfc1}) now reads:
 {\color{black}
\bea
\label{ode1}
\frac{1}{2}\sigma^2v''(x)+\hat{f}(x, v'(x))-cv(x)=0, \q x\ge 0;\qq
v(0)=0,
\eea
where the function $\hat{f}$ is defined by 
\bea
\label{fdxx1}
\hat{f}(x, z):=\Big\{ \mu(x) z +\lambda \ln\Big[\frac{\lambda(e^{\frac{a}{\lambda}(1-z)}-1 )}{1-z}\Big]\Big\}\1_{\{z\neq1\}}+[
\mu(x)+\lambda\ln a]\1_{\{z=1\}}.
\eea}
Such an ODE is only slightly more complicated than the original ODE (\ref{ode}), with the same degree of difficulty due to the transcendental nature of the equation; hence, it cannot be solved explicitly.
%, nor use the standard simulator as it would not work. 

Following the recent works \cite{MWZ1, MWZZ26}, we shall now embed our problem into an infinite horizon entropy-regularized stochastic control problem and design an RL algorithm utilizing neural networks. 
%It turns out that our framework can be imbedded into those of \cite{MWZ. MWZZ26}, as a special infinite horizon case, which we now explain.
We begin by unifying our notations with those of \cite{MWZ} and \cite{MWZZ26}.

% To compare with the PIA and NN algorithm with 
{\color{black}  First, for a given policy $\ppi\in\sA^s_{cl}$, recall the function $b^\ppi(\cd)$ and $r^\ppi(\cd)=r^\ppi_\l(\cd)$ 
%denote the {\it Shannon's entropy} by $\cH(\pi):=- \int_0^a  \pi(w) \ln \pi(w) dw$ and, in light of  
defined by (\ref{ABr}) (we often dropp  $\l$ for simplicity), and rewrite 
%, define for any $\pi\in\sP([0,a])$, $\tilde  \f( x,  \pi) := \int_0^a \f( x, w) \pi(w) dw $, for $\f(t,w)=b(x,w)=(\m(x)-w)$ and $r(x,w)=w$, respectively. Then, comparing to the functionals defined in (\ref{Jxpi}) and  (\ref{ABr}), we see that for any $\ppi\in\sA^s_{cl}$, we can write \bea
%r_\l^\ppi(t)=r^{\ppi}_\l(X^\ppi_t ):= \tilde{r}(X^\ppi_t, \ppi)-\l \cH(\ppi(\cd, X^\ppi_t)), 
%\eea
 the entropy-regularized control problem (\ref{sde2}) and (\ref{v}) as:
\bea
\label{tXpiV} 
\left\{\ba{lll}
dX^\ppi_t =  b^\ppi(X^\ppi_t) dt +\si dW_t, \q  t \ge 0; ~ X^\ppi_0 = x, \ss\\
\dis V(x)=\sup_{\ppi\in \sA^s_{cl}}\hE_x\Big[\int_0^{\t^{\ppi}}e^{-ct} r^{\ppi}_\l(X^\ppi_t ) dt\Big]. 
\ea\right.
\eea
Clearly, this is similar to the infinite horizon case in \cite{MWZ}, except that our problem is stopped at the ruin time $\t^\ppi=\t^\ppi_x$. 
%:=\inf\{t>0: X^{\ppi,x}_t<0\}$, which is in general non-differentiable with respect to the initial state $x$. 
Furthermore, if we define
%
%Now let us consider only the feedback controls: $\pi_t(w)=\pmb\pi(w, X^{\pmb\pi}_t)$, where $\pmb\pi(w,x)$ is a deterministic function, and define, as in \cite{MWZ}, the Hamiltonian by
%\bea
%\label{H}
%H(x,z)&=&\neg\sup_{{\ppi}\in\cA}\int_0^a\neg\Big[w\neg-\neg\lambda\ln\ppi(w,x)+\neg(\mu(x)-w)z\Big]\ppi(w,x)dw\nonumber\\
%&=&\neg\sup_{{\ppi}\in\cA}\int_0^a\neg\big[b(w,x)z+i(w)-\l \cH(\ppi(w,x))\big]\ppi(w,x)dw\\
%&=&\neg\sup_{{\ppi}\in\cA} \big[\tilde b(x, \ppi)z+\tilde r(x,\ppi)-\l\cH(\ppi(\cd, x))\big]. \nonumber
%\eea
%Also, we define
\bea
\label{gamma}
\g(x, z,w)=\exp\Big[\frac1{\l}[(\m(x)-w)z+w\big]\Big], 
%\ss\nonumber\\ 
%&=&e^{\frac1{\l}\m(x)z}\cd e^{\frac{w}{\l}(1-z)}, 
\q z\neq 1,
\eea
then one can easily check that the Gibbs function (see (\ref{zykz})) and function $\hat{f}$ in \reff{fdxx1} can be written as, respectively,
%the for $z\neq 1$ we can easily check that
%\bea
%\label{intg}
%\int_0^a\g(x,z,w)dw=e^{\frac1{\l}\m(x)z}\frac{\l}{(1-z)}\big[e^{\frac{a}{\l}(1-z)}-1\big].
%\eea
% Thus, recalling (\ref{Gibbs}) we derive from (\ref{gamma}) and (\ref{intg}) that
\bea
\label{Gamma}
\left\{\ba{lll}
\dis G(w,1-z)=\frac{\g(x,z,w)}{\int_0^a\g(x,z,w)dw} \1_{\{z\neq 1\}}+ \frac1{a} \1_{\{z=1\}},  \ms\\
%\eea
%%Clearly, one can easily show that $\lim_{z\to 1}\G(w, z)=$. (So we do not need to write (\ref{Gibbs}) in the ``piecewise" form!). Furthermore, as was shown in \cite{MWZ} that, 
%%Moreover,  under the optimal feedback control $\ppi^*(x,w):=G(w,1-v'(x))$, it holds that
%and 
%\bea
%\label{Hexp}
\dis \hat{f}(x,z)=\l \ln\Big[\int_0^a \g(x,z,w)dw\Big]\1_{\{z\neq 1\}}+[\mu(x)+\lambda\ln a]\1_{\{z=1\}}.
\ea\right.
\eea
%where $f(x,z)$ is defined by (\ref{fdxx}). 
%and $v$ satisfies the HJB equation:
%\bea
%\label{ode1}
%cv(x)=\frac{1}{2}\sigma^2v''(x)+\hat{f}(x, v'(x)), \q x\ge 0;\qq
%v(0)=0.
%\eea
We note that in this case the HJB equation \reff{ode1} is the same as the infinite horizon HJB equation in \cite{MWZ}, except that it is defined only on $[0,\infty)$.
%, whereas in \cite{MWZ} the HJB equation is defined on $(-\infty, \infty)$.

%\subsection

\ms 
{\bf The Bismut-Elworthy-Li Representation of $v'(x)$. }~
To explain the main idea of our algorithm, we first note that the PIA algorithm discussed in the previous sections only gives a theoretical map for an approximating sequence of the value function $V$ (or the solution  $v$ of the HJB equation (\ref{ode1})). However, each updated policy $\pi_{n+1}$ in scheme \eqref{pin} involves the approximation $v'_n$.
% In order to define the updated policy $\pi_{n+1}$ via the scheme \eqref{pin}, we need to approximate $v'_n$. 
% of the policy $\pi_n$ 
% in order to define the updated policy $\pi_{n+1}(x,w)=G(w, 1-v_n'(x))$.
A conceivable attempt would then be to 
%follow the theoretical PIA scheme and 
% first approximate the value function $v_n$, by a 
simply take the derivative of the neural network approximation, say $\hat{v}_n'$.
But such a scheme would generally fail, since even a ``good" approximation $\hat{v}_n$ based on Martingale Loss by no means implies that $\hat{v}_n'$ 
% $\{\hat{v}_n'\}_{n \geq 0}$
is a good approximation of $v_n'$ 
(see \cite{MWZZ26} for 
an example). To overcome this difficulty, we shall follow the recent works  \cite{MWZ, MWZZ26} and use an independent neural network to estimate   $v'_n$ simultaneously via the Bismut-Elworthy-Li (BEL) representation formula. The main obstacle here is that our problem is stopped at the ruin time $\t_\ppi$, which is in general non-differentiable in $x$. But this can be dealt with by a simple theoretical argument in the Markovian framework, which we now describe. 

To begin with,  for each $\ppi\in\sA^s_{cl}$ let us write
\bea
\label{J}
J(x, \ppi) &:=&  \hE_x\Big[\int_0^{\t^{\pmb\pi}} e^{-c t}  r^\ppi(X^{\pmb\pi}_t)dt\Big] = \hE_x\Big[\int_0^\infty e^{-c t}\1_{\{\t^{\pmb\pi}>t\}}  r^\ppi(X^{\pmb\pi}_t) dt\Big]. 
\eea
Since $\t^{\pmb\pi}=\t^\ppi_x$
%:=\inf\{t>0: X^{\ppi, x}_t<0\}$ 
is an $\{\cF^{X^{\pmb\pi}}_t\}$-stopping time, $\{\t^{\pmb\pi}>t\}\in \cF^{X^{\pmb\pi}}_t$, $t\in[0,T]$. By the Doob-Dynkin Lemma and then the Markov property of $X^\ppi$ we can write $ \1_{\{\t^{\pmb\pi}>t\}}=
%\Psi(X^{\pmb\pi}_{\cd\wedge t})$ where $\Psi:\hC([0,T])\mapsto \hR_+$ is some measurable functional. The Markov property  then further implies that $ \1_{\{\t^{\pmb\pi}>t\}}=
\Phi(X^{\pmb\pi}_t)$, for some measurable function $\Phi:\hR\mapsto \hR$.
%where the last equality is due to the Markov property of $X^{\pmb\pi}$. 
Thus, if $\pmb\pi^*(\cd, \cd)\in\sA^s_{cl}$ is an optimal strategy, then by (\ref{J}) we have
%\bea
%\label{J1}
%J(x, {\pmb\pi})= \hE_x\Big[\int_0^\infty e^{-c t}\Phi(X^{\pmb\pi}_t)   r^\ppi(X^{\pmb\pi}_t)dt\Big].
%\eea
%In particular,, 
%%and $\pi^*_t(w)=\pmb\pi^*(w,X^{\pmb\pi^*}_t)$, 
%then 
\bea
\label{V1}
V(x)=J(x, {\pmb\pi}^*)= \hE_x\Big[\int_0^\infty e^{-c t}\Phi(X^{\pmb\pi^*}_t) r^\ppi(X^{{\pmb\pi}^*}_t)dt\Big].
\eea
Now,  applying the Bismut-Elworthy-Li (BEL) formula \cite{Bismut, EL} or the representation formula in \cite{MaZhang1}, we obtain the probabilistic representation of $V'$ (keep  in mind that $\Phi(X^{\pmb\pi^*}_t)=\1_{\{\t^{\ppi^*}>t\}}$):
\bea
\label{V'}
 V'(x) =\hE_x\Big[\int_0^\infty \Phi( X^{\pmb\pi^*}_t)r^\ppi(X^{{\pmb\pi}^*}_t) N^{\pmb\pi^*}_t dt\Big]=\hE_x\Big[\int_0^{\t^{\ppi^*}}r^\ppi(X^{{\pmb\pi}^*}_t)N^{\pmb\pi^*}_t dt\Big].
 \eea
 where for $\ppi\in\cA$, the representation kernel $N^\ppi$ is defined by 
 {{\color{black}
 \bea
\label{utxnrep}
%\dis \pa_x u^{n}(t, x) = \dbE\Big[\pa_x g(X^{t,x}_T)\td X^{t,x}_T+\int_t^T f^n(s, X^{t,x}_s) N^{t,x}_s ds\Big];\\
\dis N^{\ppi}_t:=\frac{1}{t}\int_0^t (\sigma^{-1} \partial_x X^{\ppi}_s) dW_s; \q \partial_x  X^{\ppi}_t = 1 +  \int_0^t \pa_x b^\ppi(X^{\ppi}_s) \partial_x X^{\ppi}_s ds, \q t\in[0, T]. 
\eea}
%and we have the easy estimates
%\bea
%\label{Ntxest}
% \hE_x\big[|\td X^{\ppi}_t|^2\big] \le Ce^{C t},\q  \hE_x\big[|N^{\ppi}_t|^2\big]\le {C\over t}e^{C t}, \qq t\in(0,T].
%\eea
%
%We should note that the process $N^\ppi$ has a singularity at $t=0$ (at the order $\sim\frac{W_t}{t}$). Therefore, one cannot simple mindedly claim  from (\ref{V'}) and the fact that $\t_{\ppi}\to 0$ as $x\to 0$ that $V'(0)=0$. It should still follow the result of the paper to decide $V'(0+)=\a$, for some $\a >0$.
%it is a little subtle to analyze the behavior of $V'(x)
Our numerical algorithm will be based on two independent neural networks via the representations 
%that will  simultaneously estimate the functions $v$ and $v'$, via (
(\ref{V1}) and (\ref{V'}), respectively. 
}

\newenvironment{breakablealgorithm}
  { \begin{center}
     \refstepcounter{algorithm}%
     \hrule height .8pt depth0pt \kern 2pt
     \renewcommand{\caption}[2][\relax]{%
       {\raggedright\textbf{Algorithm \thealgorithm} ##2\par}%
       \ifx\relax##1\relax
         \addcontentsline{loa}{algorithm}{\protect\numberline{\thealgorithm}{##2}}%
       \else
         \addcontentsline{loa}{algorithm}{\protect\numberline{\thealgorithm}{##1}}%
       \fi
       \kern2pt\hrule\kern2pt
     }%
  }
  { \kern2pt\hrule\end{center} }
\makeatother

\subsection{NN Algorithm}
%\setcounter{equation}{0}

%{\color{blue} Gaozhan}: could you provide a pseudo-code for the algorithm? Here are some ideas we discussed before:
%\begin{itemize}
%\item Stop the state process $X^\ppi$ when it is under a certain threshold $\e>0$
%
%\item Estimate the derivative $V'(0+)\sim \frac{V(x)}{x}$ when $x\sim 0$.
%{\color{red}no need to use this anymore}
%\item Is it possible to use the representation formula?
%{\color{red} yes}
%\end{itemize}
% =======================
% Pseudo-code: PIA with BEL derivative (Policy Evaluation + Policy Improvement only)
% =======================

\ms

Let us now present our RL algorithm in the spirit of PIA, but involving two independent neural networks to approximate the value function and its derivative simultaneously, {\color{black}
 based on the representations (\ref{V1}) and (\ref{V'}).
%, respectively,
We denote the neural networks for $v$ and $v'$ by  $V^{\psi}$ and $P^{\eta}$  respectively, where $\psi$ and $\eta$   denote the parameters of the corresponding neural networks. }
%The detailed  In, we present
 
 \begin{rem}
 We remark that all the test results in this section were obtained on a local machine with an NVIDIA GeForce RTX 3080 Ti GPU (12 GB GDDR6X, 384-bit memory interface). For all the detailed neural network specifications, hyperparameters, and the complete codebase, please refer to our GitHub repository https://github.com/GaozhanWang/BGMW26.
 \end{rem}
 
\ss
\begin{breakablealgorithm}
\caption{Policy Iteration with Derivative Estimate (PE + PI)}
\label{alg:pia}
\begin{algorithmic}[1]\small
\State \textbf{Inputs:}
\Statex \quad Discount factor $c>0$, time step $\Delta t$, horizon cap $T$, action set $A=[0,a]$,
\Statex \quad Drift $\mu(\cdot)$, known volatility $\sigma$, and running reward $r(x,w)$,
\Statex \quad Temperature parameter $\lambda>0$,
\Statex \quad Number of initial wealth states $B$, Monte Carlo paths per state $M$, and optimizer minibatch size $\mathcal B$,
\Statex \quad PE training epochs $E_V,E_P$, and learning rates
$\eta_V,\eta_P$,
\Statex \quad Neural networks $V^{\psi}(x)$ and $P^{\eta}(x)$,
\Statex \quad Maximum PI iterations $N_{\mathrm{PI}}$, patience
$P_{\mathrm{PI}}$, and tolerance
$\varepsilon_{\mathrm{PI}}$,
\Statex \quad State grid $\mathcal G$
\State \textbf{Output:} Learned value $V^{\psi}$, learned derivative
$P^{\eta}$, and improved policy $\pi(\cdot\mid x)$.

\State \textbf{Initialize:} Randomly initialize a Gibbs policy
$\pi^{(0)}(w\mid x)=G(w,\lambda\beta^{(0)}(x))$, and  $(\psi^{(0)},\eta^{(0)})$; set $n\gets0$ and $q\gets0$.

\While{$n<N_{\mathrm{PI}}$ \textbf{and} $q<P_{\mathrm{PI}}$}

  \State \textbf{(PE-1) Generate training trajectories:}
\State Hold the current policy $\pi^{(n)}$ fixed.
\State Uniformly sample $B$ distinct states
$\{x_0^{(b)}\}_{b=1}^{B}\subset\mathcal G$ without replacement.
\State Set $N_{\mathrm{train}}\gets BM$.

\For{$b=1,\dots,B$}
  \For{$j=1,\dots,M$}
    \State Set $m\gets(b-1)M+j$,
    $X_0^{(m)}\gets x_0^{(b)}$, and $k\gets0$.

    \While{$k<N_T$ \textbf{and} $X_{t_k}^{(m)}>0$}
      \State Sample
      $w_k^{(m)}\sim
      \pi^{(n)}(\cdot\mid X_{t_k}^{(m)})$.
      \State Independently draw
      $\Delta W_k^{(m)}\sim\mathcal N(0,\Delta t)$.
      \State Euler step:
      $$
      X_{t_{k+1}}^{(m)}
      \gets
      X_{t_k}^{(m)}
      +
      \bigl(\mu(X_{t_k}^{(m)})-w_k^{(m)}\bigr)\Delta t
      +
      \sigma\Delta W_k^{(m)}.
      $$
      \State $k\gets k+1$.
    \EndWhile

    \State Set $K^{(m)}\gets k$.
  \EndFor
\EndFor

  \State \textbf{(PE-2) Build value targets:}
  \For{$m=1,\dots,N_{\mathrm{train}}$}
    \State
    \[
      \widehat{\mathcal U}^{(m)}
      \gets
      \sum_{k=0}^{K^{(m)}-1}
      e^{-ct_k}
      r\bigl(X^{(m)}_{t_k},w^{(m)}_k\bigr)\Delta t
    \]
  \EndFor

  \State \textbf{(PE-3) Train the value network $V^{\psi}$:}
  \For{$\mathrm{epoch}=1,\dots,E_V$}
    \For{each minibatch
$\mathcal B\subset\{1,\dots,N_{\mathrm{train}}\}$
of size $\mathcal B$}
      \State
      \[
        L_V
        \gets
        \frac{1}{|\mathcal B|}
        \sum_{m\in\mathcal B}
        \left|
          V^{\psi}(X^{(m)}_0)-\widehat{\mathcal U}^{(m)}
        \right|^2
      \]
      \State
      $\psi\gets\psi-\eta_V\nabla_{\psi}L_V$
    \EndFor
  \EndFor

  \State \textbf{(PE-4) Build derivative targets:}
  \For{$m=1,\dots,N_{\mathrm{train}}$}
    \State Compute the kernel process
    $\{N^{(m)}_{t_k}\}_{k=0}^{K^{(m)}}$
    from the derivative representation formula
    \State
    \[
      \widehat{\mathcal P}^{(m)}
      \gets
      \sum_{k=0}^{K^{(m)}-1}
      e^{-ct_k}
      r\bigl(X^{(m)}_{t_k},w^{(m)}_k\bigr)
      N^{(m)}_{t_k}\Delta t
    \]
  \EndFor

  \State \textbf{(PE-5) Train the derivative network $P^{\eta}$:}
  \For{$\mathrm{epoch}=1,\dots,E_P$}
   \For{each minibatch
$\mathcal B\subset\{1,\dots,N_{\mathrm{train}}\}$
of size $\mathcal B$}
      \State
      \[
        L_P
        \gets
        \frac{1}{|\mathcal B|}
        \sum_{m\in\mathcal B}
        \left|
          P^{\eta}(X^{(m)}_0)-\widehat{\mathcal P}^{(m)}
        \right|^2
      \]
      \State
      $\eta\gets\eta-\eta_P\nabla_{\eta}L_P$
    \EndFor
  \EndFor

  \State \textbf{(PI) Policy improvement via the Gibbs form:}
  \State Define
  \[
    \beta^{(n+1)}(x)
    \gets
    \frac{1-P^{\eta}(x)}{\lambda},
    \qquad
    \pi^{(n+1)}(w\mid x)
    \gets
    G\bigl(w,1-P^{\eta}(x)\bigr).
  \]

  \State Compute the maximum policy-parameter change
  \[
    d_n
    \gets
    \max_{x\in\mathcal G}
    \left|
      \beta^{(n+1)}(x)-\beta^{(n)}(x)
    \right|.
  \]

  \If{$d_n>\varepsilon_{\mathrm{PI}}$}
    \State $q\gets0$
  \Else
    \State $q\gets q+1$
  \EndIf

  \State $n\gets n+1$

\EndWhile

\State \textbf{Final PE:} With $\pi^{(n)}$ fixed, perform one final
policy-evaluation step so that $V^{\psi}$ and $P^{\eta}$ correspond to
the returned policy.

\State \textbf{Return}
$\bigl(V^{\psi},P^{\eta},\pi^{(n)}\bigr)$.
\end{algorithmic}
\end{breakablealgorithm}

{\color{black}We should remark that in Algorithm 1, instead of 
%In the Gibbs-form policy update, an accurate approximation of the spatial is critical. Rather than 
calling for the derivative of the value function from the neural network output $V^\psi$, 
%which can be essentially wrong because policy evaluation targets the value itself, not its derivative, 
we estimate the derivative directly using a separate neural network $P^\eta$, based on the  BEL formula. 

Specifically, both the value function approximation $V^{\psi}$  and its derivative $P^{\eta}$ are represented by  independent, fully connected feedforward neural networks, with the following main features:
\vspace{-2mm}
\begin{itemize}
% \item{} Each network uses four hidden layers with 128 neurons per layer, and the hyperbolic tangent activation function is applied at all hidden layers.
\item We use an equally spaced state grid
$\mathcal G=\{x_j\}_{j=1}^{400}\subset[0.05,3]$.
At each policy-evaluation step, the current policy $\pi^{(n)}$
is held fixed. We uniformly sample $B=256$ distinct initial
wealth states from $\mathcal G$ without replacement.
For each selected state, we simulate $M=128$ independent
trajectories under $\pi^{(n)}$, using independent action
randomizations and Brownian increments. Thus, each
policy-evaluation step generates $BM=32{,}768$ training
trajectories.

Each trajectory provides one value target and one derivative
target, paired with its initial wealth. These targets are
constructed once and reused throughout the corresponding
policy-evaluation step. The value network $V^\psi$ and the
derivative network $P^\eta$ are trained separately using Adam,
for at most $E_V=E_P=200$ epochs, with learning rates
$\eta_V=\eta_P=8\times10^{-4}$.
% \item We consider an equally spaced  State grid $\mathcal G=\{x_j\}_{j=1}^{400}\subset[0.05,3]$. We randomly initialize $x_0 \in \mathcal{G}$, and simulate $M=128$
% Monte Carlo trajectories. At each policy evaluation, we use $B=256$ initial wealth states, and for each state we simulate $M=128$
% Monte Carlo trajectories. Thus, a total of $BM=32{,}768$ trajectories
% are generated at each epoch to construct the value and derivative targets.
% The value network $V^{\psi}$ and the derivative network $P^{\eta}$
% are trained simultaneously using the Adam optimizer for at most
% $ E_V=E_P=200$ epochs, with learning rates $ \eta_V=\eta_P=8\times10^{-4}$.

% More precisely, at each policy-evaluation step, validation is performed on $256$
% equally spaced wealth states in $\mathcal G=\{x_j\}_{j=1}^{400}\subset[0.05,3]$, using $128$ independently
% simulated trajectories per state. 

\item Training of either network is terminated
early if its validation loss fails to improve by at least $10^{-4}$ for
$12$ consecutive epochs, and the network parameters attaining the
smallest validation loss are retained. The outer policy-iteration loop is terminated either after
($N_{\mathrm{PI}}=50$) policy-improvement steps or when the maximum
policy-parameter change does not exceed
($\varepsilon_{\mathrm{PI}}=10^{-4}$) for
($P_{\mathrm{PI}}=12$) consecutive iterations, whichever occurs first.
Upon termination, one final policy-evaluation step is performed with the
learned policy held fixed, so that the reported $V^{\psi}$ and
$P^{\eta}$ correspond to the final policy. 

\vspace{-2mm}
\item Although the underlying dividend problem has an infinite time horizon, in the algorithm we follow the strategy of ``batch learning" and truncate the simulated trajectories at a finite terminal time $T>0$ (which we set as $T=2$). Such a choice of $T$ turns out to be sufficient due to the relatively large discount factor $c=10$, which effectively suppresses long-term contribution when we generate multiple paths for the aggregated results.   (In light of Remark \ref{remark5}, given $c=10$ and $\D t =0.02$, $T \geq 0.2$ would already be sufficient.)

\vspace{-2mm}
\item To explicitly monitor the ruin time, 
%defined as the first hitting time at which the state process falls below $0$, 
for each simulated path, we set the ``ruin" threshold $\e=10^{-3}$ if it happens before the truncation time $T$. In other words, 
%ed accordingly, while for those that do not, the remaining 
the tail contribution beyond $\t^{(m)}\wedge T$, where $\t^{(m)}$ is the approximated ruin time at each iteration, is neglected.
%, resulting in a controlled truncation error.

\end{itemize}

{\color{black} We should note that, although the trajectories in our algorithm are simulated
with given $\m$ and $\si$, the training uses only the sampled trajectory values, and not
the coefficients themselves. Consequently, the algorithm remains applicable when only the
trajectories are observed, and the model parameters are unknown.
}
We would like to emphasize that the Algorithm 1 proposed above is applicable to a fairly large class of diffusion models beyond the constant-coefficient case. Some general convergence results have been established in recent works \cite{MWZ, MWZZ26}, once we embed our model into the general infinite horizon case, as we outlined at the beginning of this section.

\subsection{Numerical Results}

We now present two numerical examples: one for constant coefficients and one for state-dependent coefficients. 
%The results are equally effective. 
We should note that since the entropy-regularized HJB equation (\ref{ode1})  is a
%. This method and the reasoning are discussed in \cite{MWZZ26}. We want to claim that due to the 
transcendental equation, it is not possible to find the explicit solution. Therefore, we cannot assess the convergence by the usual error (between the approximation and the true solution). Instead, we shall 
% it's not convenient to design an example with analytical true solution, we only 
look at loss decay in the PI (policy evaluation) step and track the value function for the policy improvement.  As we can see, the result is quite satisfactory.
}
\begin{eg}\label{EX1}
In this example, we consider the simple constant coefficient case. We shall set the  system  parameters $\si=1$ and 
$\mu(x) \equiv 3$. Also, we set $x_0=1$, $c=10$, $T=2$, $\Delta t=0.02$.  
\begin{figure}[H]
    \centering
    \begin{subfigure}[t]{0.5\textwidth}
        \centering
        \includegraphics[width=\textwidth]{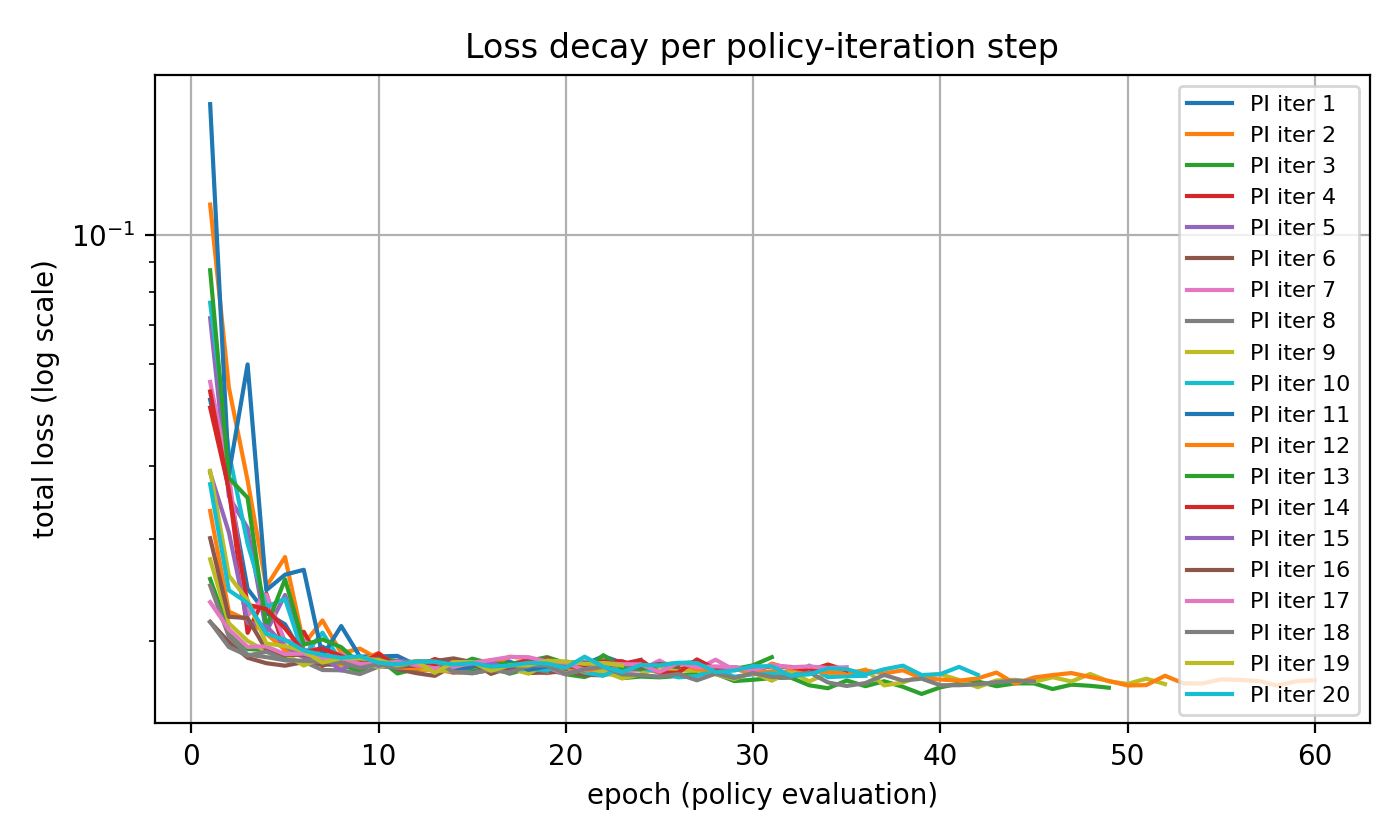}
        \caption{Loss decay during policy evaluation.}
        \label{fig:pe1}
    \end{subfigure}
    \hskip 12pt
    \begin{subfigure}[t]{0.4\textwidth}
        \centering
        \includegraphics[width=\textwidth]{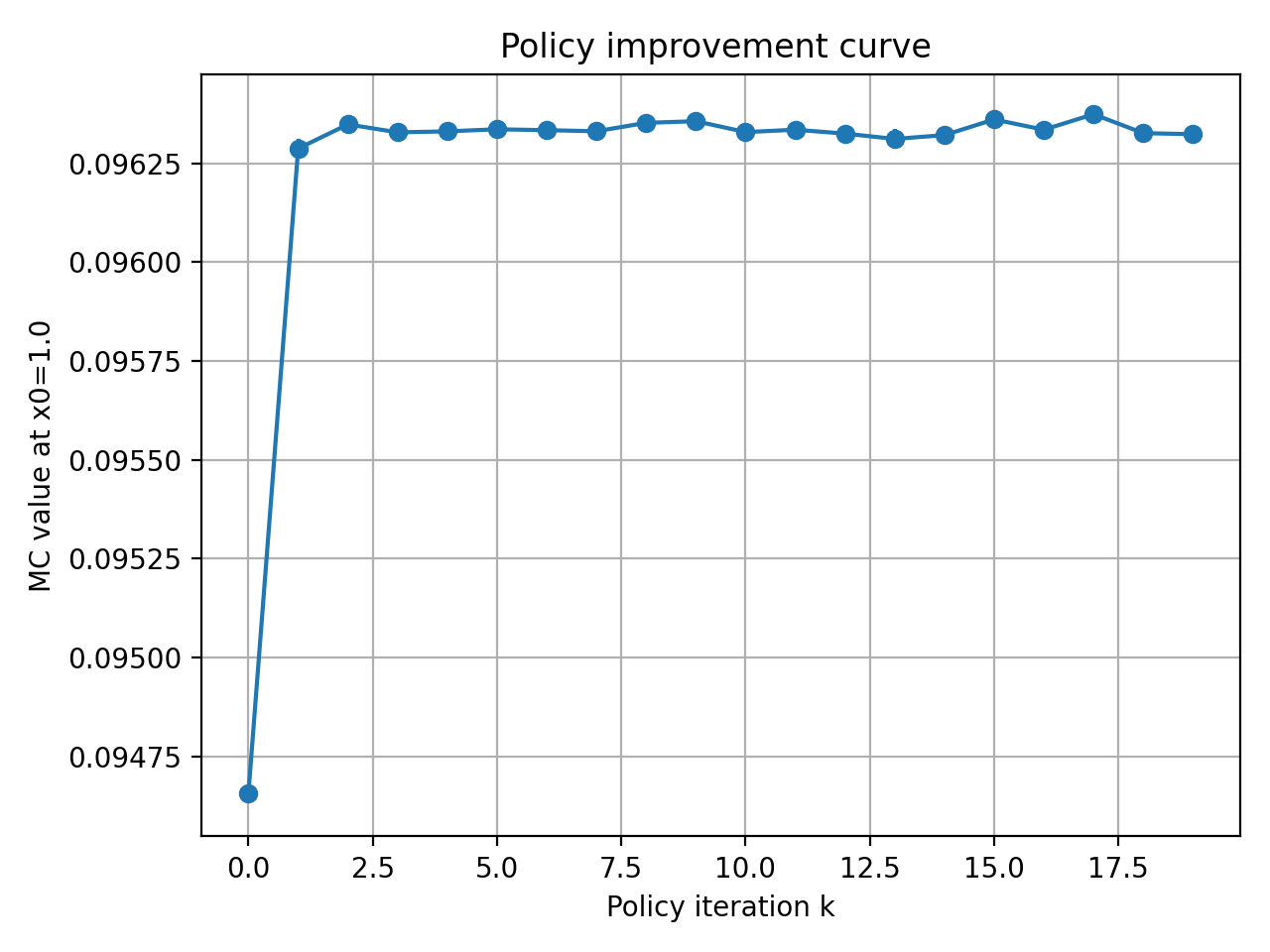}
        \caption{\footnotesize Value improvement along the policy-iteration sequence.}
        \label{fig:pi1}
    \end{subfigure}

    \caption{\small Convergence behavior of policy evaluation and policy improvement
    in Example~\ref{EX1}.}
    \label{fig:pe-pi-ex1}
\end{figure}
\vspace{-2mm}
{\rm Figure~\ref{fig:pe-pi-ex1} summarizes the numerical convergence behavior of
the PIA scheme. Panel~(\subref{fig:pe1}) shows the representative
decay of the training losses during the policy-evaluation steps.
Panel~(\subref{fig:pi1}) reports the estimated value function at the fixed
evaluation point \(x_0\) along the policy-iteration sequence. {\color{black} Its overall
increase and subsequent stabilization indicate that the successive policy
updates do improve the value objective, and the algorithm converges empirically. } 
% \begin{figure}[H]
%     \centering
%     \includegraphics[width=0.7\textwidth]{PE1.png}
%     \caption{Policy Evaluation in Example \ref{EX1}}
%     \label{fig: PE1}
% \end{figure}
% Figure \ref{fig: PE1} represents the general behavior of the loss decay in the policy evaluation steps. 
% \begin{figure}[H]
%     \centering
%     \includegraphics[width=0.55\textwidth]{PI1.png}
%     \caption{Policy Improvement in Example \ref{EX1}}
%     \label{fig: PI1}
% \end{figure}
% Figure \ref{fig: PI1} shows that the value function increases at 
% a given fixed point $x_0$, along the PI sequence, and converges.

\begin{figure}[H]
    \centering
    \includegraphics[width=0.6\textwidth]{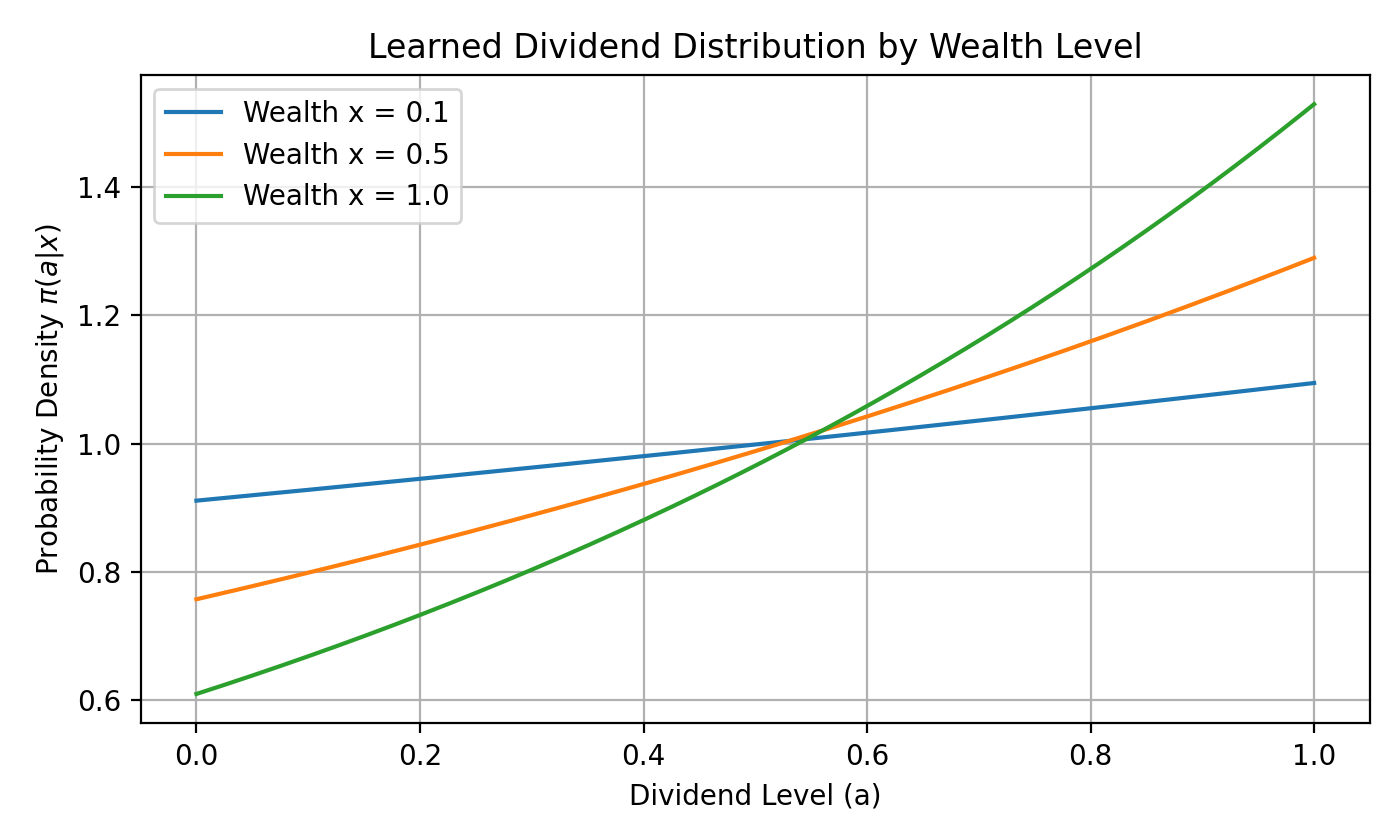}
    \caption{\small Learned optimal policy densities in Example~\ref{EX1}.}
    \label{fig: policy_pdf_wealth_const}
\end{figure}
Figure~\ref{fig: policy_pdf_wealth_const} displays the distribution of the optimal
dividend rate conditional on the current wealth levels $x\in\{0.1,0.5,1.0\}$. Each curve corresponds to one current wealth level.
\begin{figure}[H]
    \centering
    \includegraphics[width=1\textwidth]{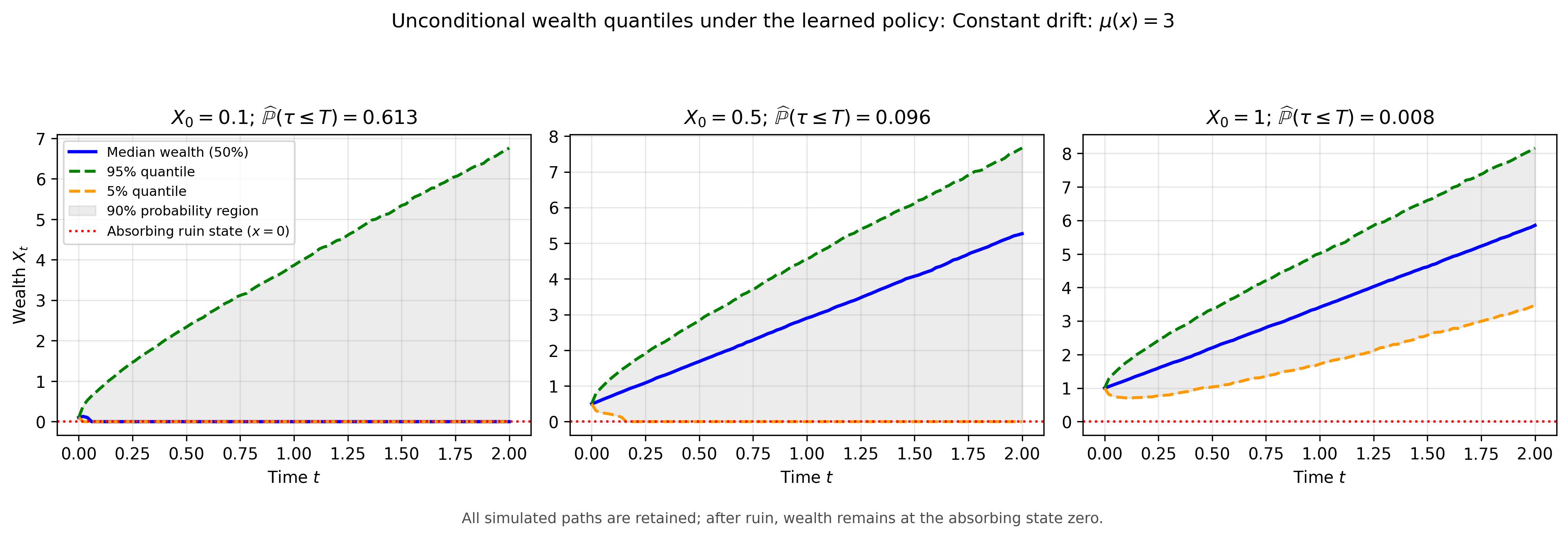}
    \caption{\footnotesize Unconditional wealth quantiles under the learned policy in
Example~\ref{EX1}, for $X_0\in\{0.1,0.5,1.0\}$.}
    \label{fig: state_quantiles_including_ruin_by_initial_wealth_const}
\end{figure}
\vspace{-2mm}
Figure~\ref{fig: state_quantiles_including_ruin_by_initial_wealth_const} describes the unconditional
distribution of the controlled wealth process. Ruined paths remain in the
sample and contribute a probability mass at zero. Consequently, the
$5\%$ quantile becomes zero once the cumulative ruin probability reaches
$5\%$, while the median becomes zero if it reaches $50\%$. This
representation therefore displays both the evolution of wealth and the effect
of ruin in a single figure.
\begin{figure}[H]
    \centering
    \includegraphics[width=1\textwidth]{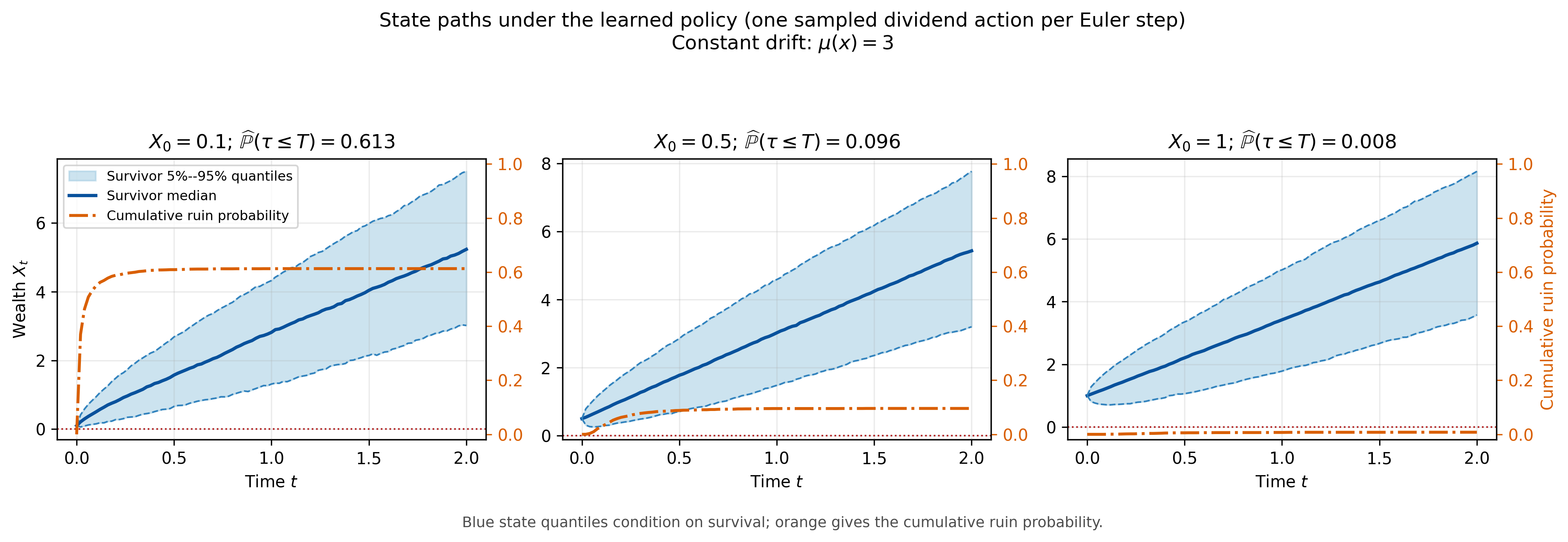}
    \caption{\footnotesize Survival-conditioned wealth quantiles and cumulative ruin
probability under the learned policy in Example~\ref{EX1}, for $X_0\in\{0.1,0.5,1.0\}$.}
    \label{fig: state_quantiles_by_initial_wealth_const}
\end{figure}
\vspace{-2mm}
Figure~\ref{fig: state_quantiles_by_initial_wealth_const} separates two different features of
the controlled process. The blue quantiles describe the wealth distribution
among paths that have survived up to time $t$, whereas the orange curve is
computed from all simulated paths and records the cumulative probability of
ruin.}
\end{eg}

We remark that Example~\ref{EX1} can be considered only as a ``sanity check" for the proposed algorithm, as it does not fully demonstrate the advantage of the RL approach. In fact, in the case of constant coefficients, standard parameter estimation methods or other existing simulators might be applied directly to approximate the solution to the original 
 control problem without using the “entropy-regularized" RL scheme at all. In the following example, we consider a setting with a state-dependent drift term, where such classical estimation techniques are no longer valid in general. This example highlights the benefit of our reinforcement learning framework in handling more general and structurally complex models. 

%{\color{red} (Comments on truncation $T$ and implement of $\t^\ppi$?)}

\ss
\begin{eg} 
\label{EX2}
{\color{black} \rm In this example, we test our algorithm in the case when the drift coefficient is state-dependent. Specifically, we set $x_0=1, \mu(x) =2x+3, c=10, T=2, \Delta t=0.02$. As we can see below, the resulting 
%Comparing the two examples we see that the 
numerical behavior of the algorithm shows that
%the two examples we see that in both cases the loss decay plots show that 
the training procedures remain stable and converge rapidly, compared to the constant coefficient case in Example  \ref{EX1}.  However, the policy improvement curve reveals a markedly faster convergence in the constant-coefficient case, due mainly to the fact that in
% the significant simplified 
the case of constant drift, the spatial derivative vanishes and does not contribute to the kernel $N$ in the BEL representation formula, 
%. In fact, in that case the policy improvement step becomes effectively close to a one-step update.
whereas in the state-dependent case, the presence of such derivative leads to 
%(Example~\ref{EX2}) 
% the derivative of the drift plays an active role, and thus the algorithm requires 
more policy-iteration steps to identify the optimal strategy. Nevertheless, the convergence of the value function is still clearly observed with a satisfactory speed.

In what follows, we present four testing results that are state-dependent counterparts to Figures 1-4. 
%Figure \ref{fig: state_quantiles_by_initial_wealth_const}. 
We believe that these results demonstrate that the proposed algorithm is both stable and robust across different model specifications.
}

\begin{figure}[H]
\centering
\begin{subfigure}[t]{0.54\textwidth}
\centering
\includegraphics[width=\textwidth]{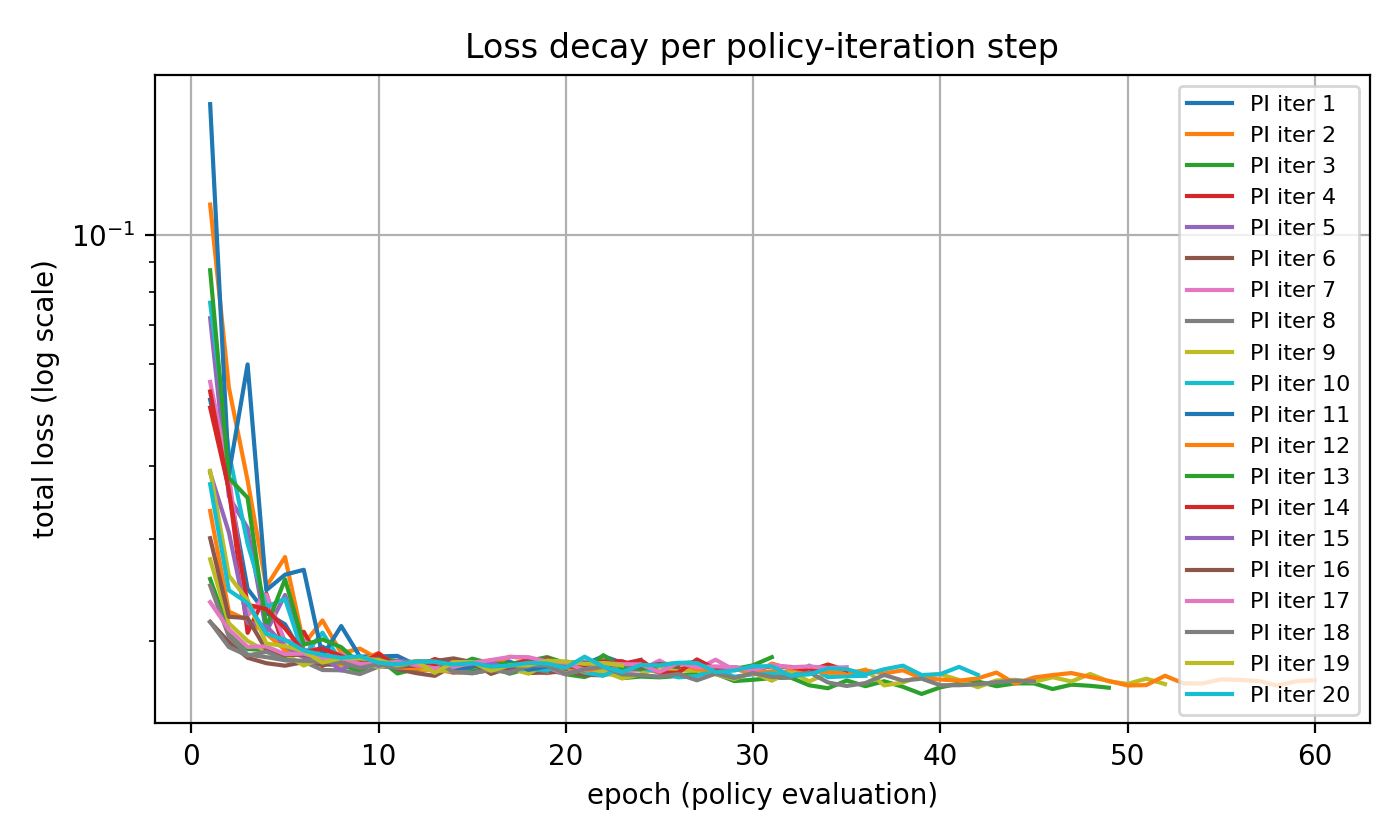}
\caption{Loss decay during policy evaluation.}
\label{fig:pe2}
\end{subfigure}
\hfill
\begin{subfigure}[t]{0.42\textwidth}
\centering
\includegraphics[width=\textwidth]{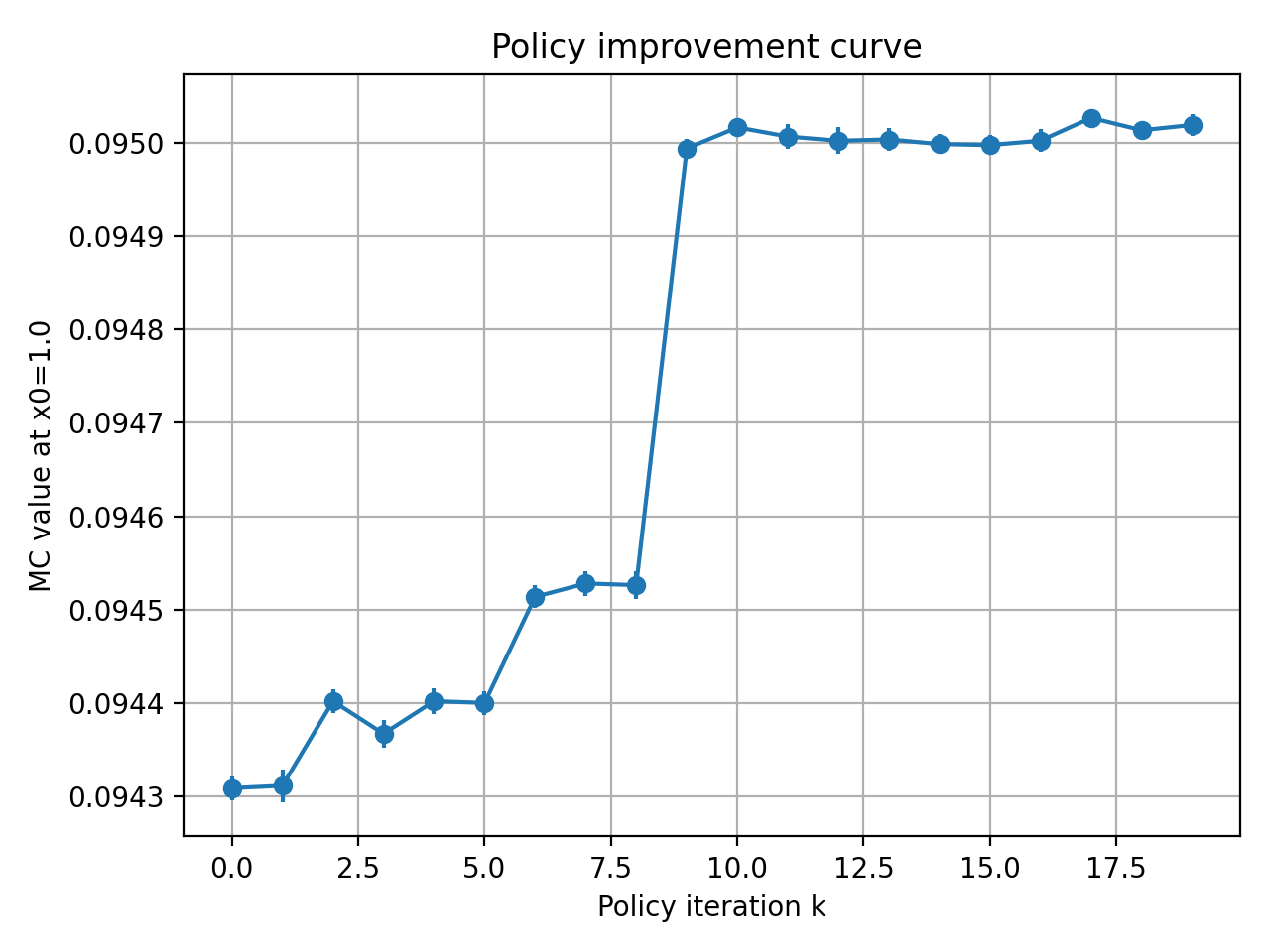}
\caption{Value improvement along the policy-iteration sequence.}
\label{fig:pi2}
\end{subfigure}
\caption{Convergence behavior of policy evaluation and policy improvement
in Example~\ref{EX2}.}
\label{fig:pe-pi-ex2}
\end{figure}
\vspace{-2mm}
%{\rm Figure~\ref{fig:pe-pi-ex2} summarizes the numerical convergence behavior of
%the policy-iteration algorithm in Example~\ref{EX2}.
%Panel~\subref{fig:pe2} shows the representative decay of the training losses
%during the policy-evaluation steps. Panel~\subref{fig:pi2} reports the
%estimated value function at the fixed evaluation point $x_0$ along the
%policy-iteration sequence. Its overall increase and subsequent stabilization
%indicate that the successive policy updates improve the objective and that the algorithm approaches convergence.}

% \begin{figure}[H]
%     \centering
%     \includegraphics[width=0.65\textwidth]{PE.png}
%     \caption{Policy Evaluation in Example \ref{EX2}}
%     \label{fig: PE2}
% \end{figure}

% \begin{figure}[H]
%     \centering
%     \includegraphics[width=0.55\textwidth]{PI.png}
%     \caption{Policy Improvement in Example \ref{EX2}}
%     \label{fig: PI2}
% \end{figure}

\begin{figure}[H]
    \centering
    \includegraphics[width=0.5\textwidth]{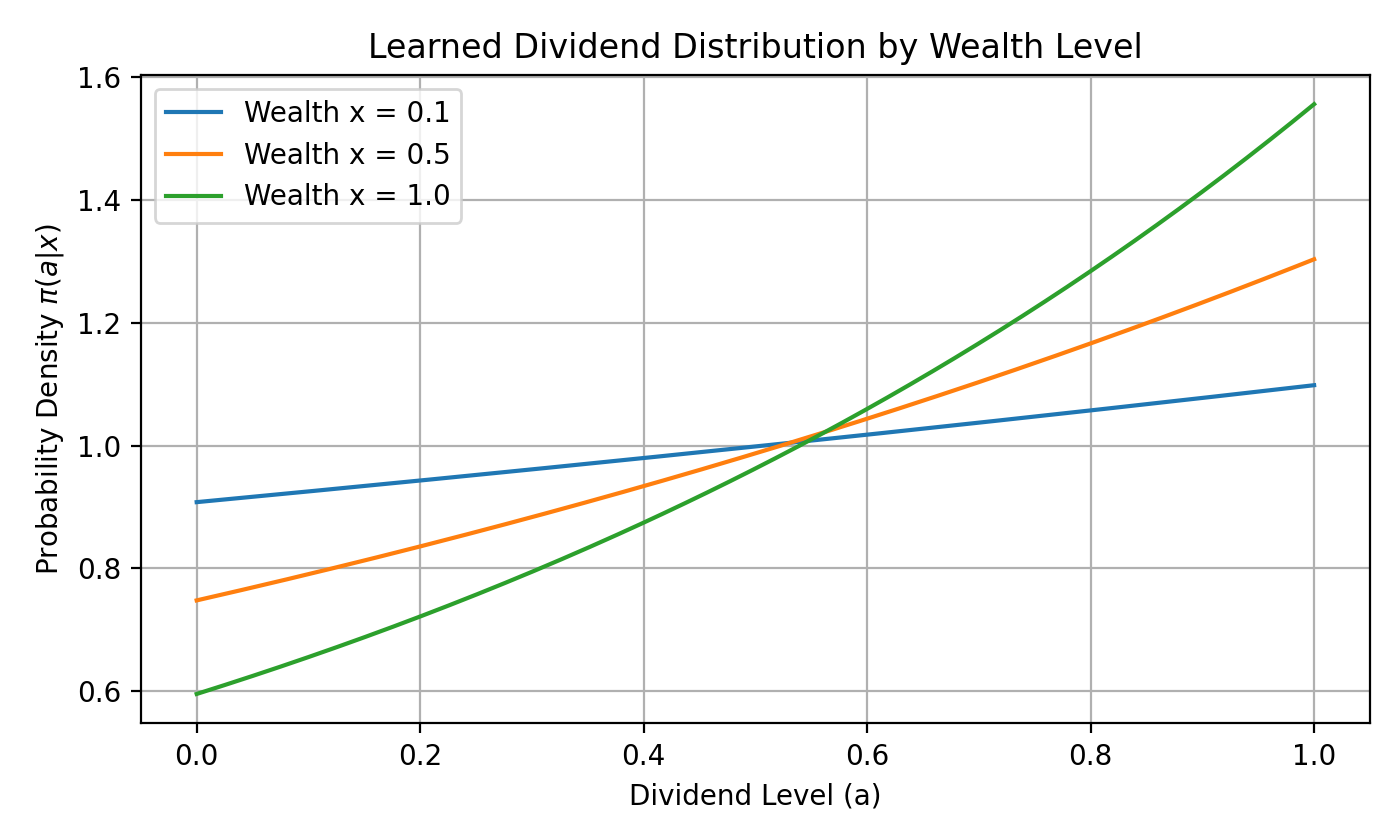}
    \caption{\small Learned optimal policy densities in Example~\ref{EX2} at the current wealth levels $x\in{0.1,0.5,1.0}$.}
    \label{fig: policy_pdf_wealth_state_dependent}
\end{figure}
\begin{figure}[H]
    \centering
    \includegraphics[width=1\textwidth]{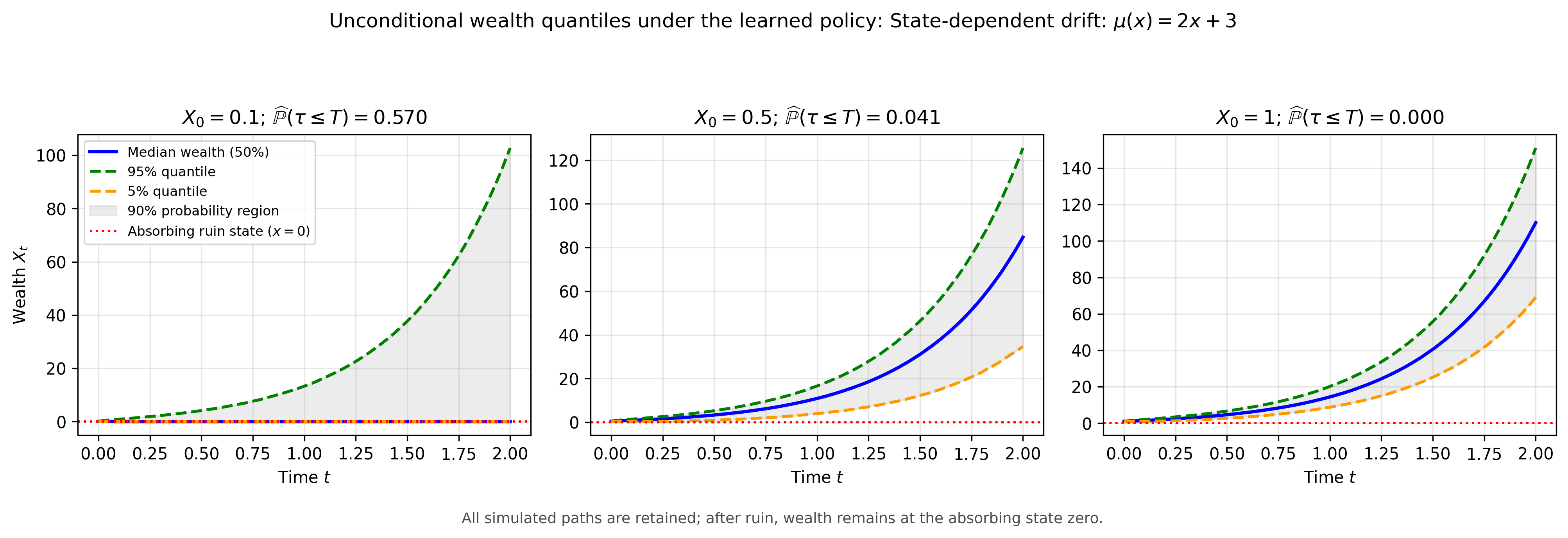 }
    \caption{\footnotesize Unconditional wealth quantiles under the learned policy in
Example~\ref{EX2}, for $X_0\in{0.1,0.5,1.0}$. The $5\%, 50\%, 95\%$ quantiles are computed over all simulated paths, including paths
absorbed at zero after ruin.}
    \label{fig: state_quantiles_including_ruin_by_initial_wealth_state_dependent}
\end{figure}
\begin{figure}[H]
    \centering
    \includegraphics[width=1\textwidth]{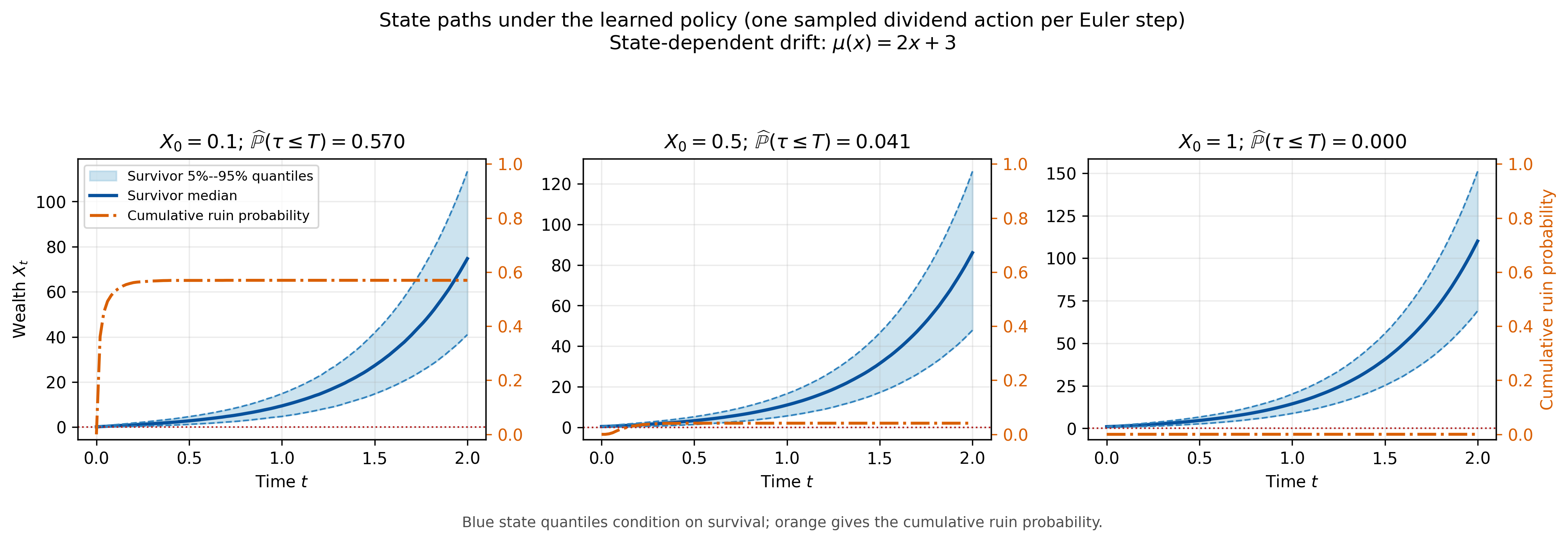}
    \caption{\footnotesize Survival-conditioned wealth quantiles and cumulative ruin
probability under the learned policy
% in Example~\ref{EX2}, 
for $X_0\in{0.1,0.5,1.0}$. Blue denotes the survivor $5\%-95\%$
quantile range and median, and orange denotes
$\widehat{\mathbb P}(\tau\leq t)$.}
    \label{fig: state_quantiles_by_initial_wealth_state_dependent}
\end{figure}
\end{eg}

}
%\begin{rem}
{\color{black} The numerical experiments in both examples have some natural financial implications that are not theoretically obvious due to the transcendental nature of the entropy-regularized HJB equation. For example, in Figures 2 and 6 we plot the optimal policy densities $\pi(\,\cd\, |x)$ with different reserve level $x$, using the neural network approximation $P^{\eta}(x)$. Both results show that when $x$ is small (near ruin), the density $\pi(\,\cd\,|x)$ becomes ``flatter", closer to the uniform distribution (the maximum-entropy distribution on a compact action set). Note that $P^\eta(x)$ approximates the marginal wealth $V'(x)$, and by Gibbs form $\pi^*(w \,|\,x) \propto \exp{\frac{(1-P^\eta(x))w}{\lambda}}$, $w\in[0,a]$, this indicates that  $(1-P^\eta(x)) \sim 0$ (or $P^\eta(x)\sim 1$). Namely, when $x$ is small, the net advantage of paying dividends is also small, and the entropy term becomes relatively more influential. As  $x$ increases, the firm has a larger solvency buffer and the marginal value of retaining additional reserves decreases, 
%. The coefficient $(1-P^\eta(x))$ therefore increases, 
the probability mass becomes more concentrated on higher dividend rates. Economically, our learned policy shown in 
Figures 2 and 6 shows a clear transition from cautious and more diversified dividend decisions near ruin to a stronger preference for immediate dividend payments at higher wealth levels.
%Specifically, since the optimal Gibbs density satisfies
%(or the derivative of value function $V'(\cd)$), 
%there is a visible wealth-dependent shift in the policy density. 
% When the firm is close to ruin, an additional unit of retained wealth is particularly valuable, and the numerical estimate  $P^\eta(x)$ is close to one. and pulls the policy density toward the uniform distribution, 
We should note that the behavior of marginal  $V'(x)$ near ruin ($x\sim 0$) does not seem to be easily obtainable from the entropy-regularized HJB equation (\ref{hjbfc1}) alone without the explicit solution  (in Proposition \ref{classic} we only argued that $V'(0)>0$ exists) 

%\begin{rem}
To end this section, we make some further comments on the independent neural network approximation for the derivative (marginal) wealth, $P^\eta$, and discuss its numerical consistency. 
%\label{rem:derivative-consistency}
As we mentioned before, the use of a separate neural network for the derivative, rather than {\it automatic differentiation} (Autodiff) of the neural network value $V^{\psi}$, is deliberate. The main reason is that the  {\it Autodiff} $\partial_xV^{\psi}$ 
%  evaluates the derivative of a fitted value network exactly, such a 
is not an accurate approximation of the derivative of the true value function.  In fact,  since the  Martingale Loss controls only the 
%approximation error of the value function in an 
$L^2$ error, which by no means implies that it can control the (Sobolev) $H^1$-error. 
%Thus, even a small value-approximation error can contain oscillations whose derivatives remain large. 
This issue would be particularly relevant in the presence of other factors such as Monte Carlo noise, a ruin
boundary, and state-dependent coefficients. 
%Consequently, convergence of the value approximation alone does not justify using in the policy-improvement step. 
Following \cite{MWZ, MWZZ26}, our algorithm uses the independent network $P^{\eta}$ and trains it directly against a stochastic
representation of the derivative obtained from the
BEL formula, so that it has its own
derivative-level training target. 
% and is not obtained merely by differentiating the value approximation. 
This is precisely the motivation discussed in \S6.1.

In general, the use of two independent networks may nevertheless create a finite-sample
discrepancy between $P^{\eta}$ and $\partial_xV^{\psi}$. In the following example, we perform
a {\it fixed-policy diagnostic} by directly examining the consistency and accuracy of the two derivative estimates.
\begin{eg}
\label{EX3}
% {\color{red}(This experiment beats the second alternative suggested by the referee)} 
{\rm We consider a simple entropy-regularized problem (\ref{v}) and (\ref{sde2}) with constant coefficients $\mu=3$, $\sigma=1$, $c=10$, $\lambda=1$, and $A=[0,1]$. Let us fix a feedback Gibbs policy
$$ \overline{\pi}(w)=     \frac{\exp(\overline{\beta}w)}
    {\displaystyle\int_0^1\exp(\overline{\beta}u)\,du}, \qq \overline{\beta}=0.75.
 $$
Denoting $\overline w = \int_0^1 w\,\overline{\pi}(w) dw$ and $\overline r =\int_0^1[w-\lambda\log\overline{\pi}(w)]\overline{\pi}(w)dw$, we can write the corresponding value function as
$V^{\overline{\pi}}(x) =  \mathbb E_x\left[ \int_0^\tau e^{-ct}\overline r\,dt\right]$. Furthermore, as in our previous experiments, we shall truncate the problem to the finite horizon with $T=2$. In light of Remark \ref{remark5}, such a choice of $T$ should make the truncation error 
negligible given the discounting factor $c=10$,

%Freezing the policy ensures that $P^{\eta}$, $\partial_xV^{\psi}$, and the
%external benchmark all correspond to the same linear fixed-policy problem.  {\color{red}see the discussion in Remark \ref{trunc}}.
%For an independent benchmark, we solve the fixed-policy boundary-value problem(BVP)
%$$
%    \frac{\sigma^2}{2}
%    \bigl(V_{\mathrm{ref}}^{\overline{\pi}}\bigr)''(x)
%    +
%    \overline b
%    \bigl(V_{\mathrm{ref}}^{\overline{\pi}}\bigr)'(x)
%    -
%    cV_{\mathrm{ref}}^{\overline{\pi}}(x)
%    +
%    \overline r
%    =0,
%    \qquad x>0,
%$$
%subject to
%$$
%    V_{\mathrm{ref}}^{\overline{\pi}}(0)=0,
%    \qquad
%    \lim_{x\to\infty}
%    V_{\mathrm{ref}}^{\overline{\pi}}(x)
%    =
%    \frac{\overline r}{c}.
%$$
Next, let us consider the corresponding HJB equation (\ref{hjbfc}) and denote its solution by $V_{\mathrm{ref}}^{\overline{\pi}}$, which is  
obtained by a high-accuracy PDE solver. We shall use this reference solution as a benchmark to compare the performance of neural network approximation $P^\eta$ and the autodiff $\partial_xV^{\psi}$. We emphasize here that the reference solution $V_{\mathrm{ref}}^{\overline{\pi}}$  will not be 
% and is used only for evaluation and is never 
supplied to learn $V^\psi$ and $P^\eta$, which are trained independently for five different random seeds.

In this experiment we provide two complementary comparisons. One compares the derivative of the reference solution $V_x^{\mathrm{ref}}$ with the independently trained derivative network $P^{\eta}$, and the other
%. The blue curve is
 %the mean of $P^{\eta}$ over five independent training seeds, and 
 %the blue shaded region is its pointwise $95\%$ confidence interval. Panel~(b) 
 compares  the across-seed mean absolute errors
$\lvert P^{\eta}-V_x^{\mathrm{ref}}\rvert$ and $\lvert\partial_xV^{\psi}-V_x^{\mathrm{ref}}\rvert$ (over five independent training seeds) on a logarithmic scale.   This example empirically validates that $P^\eta$ actually performs better than the autodiff.
% We shall see that
% $P^\eta$ actually performs better than the autodiff in this example. 
%The corresponding shaded regions are pointwise $95\%$ bootstrap confidence intervals for the mean absolute-error curves. 
We emphasize again that the reference solution $V^{\mathrm{ref}}$  is used only for evaluation.
\begin{figure}[H]
    \centering
    \includegraphics[width=0.9\textwidth]
    {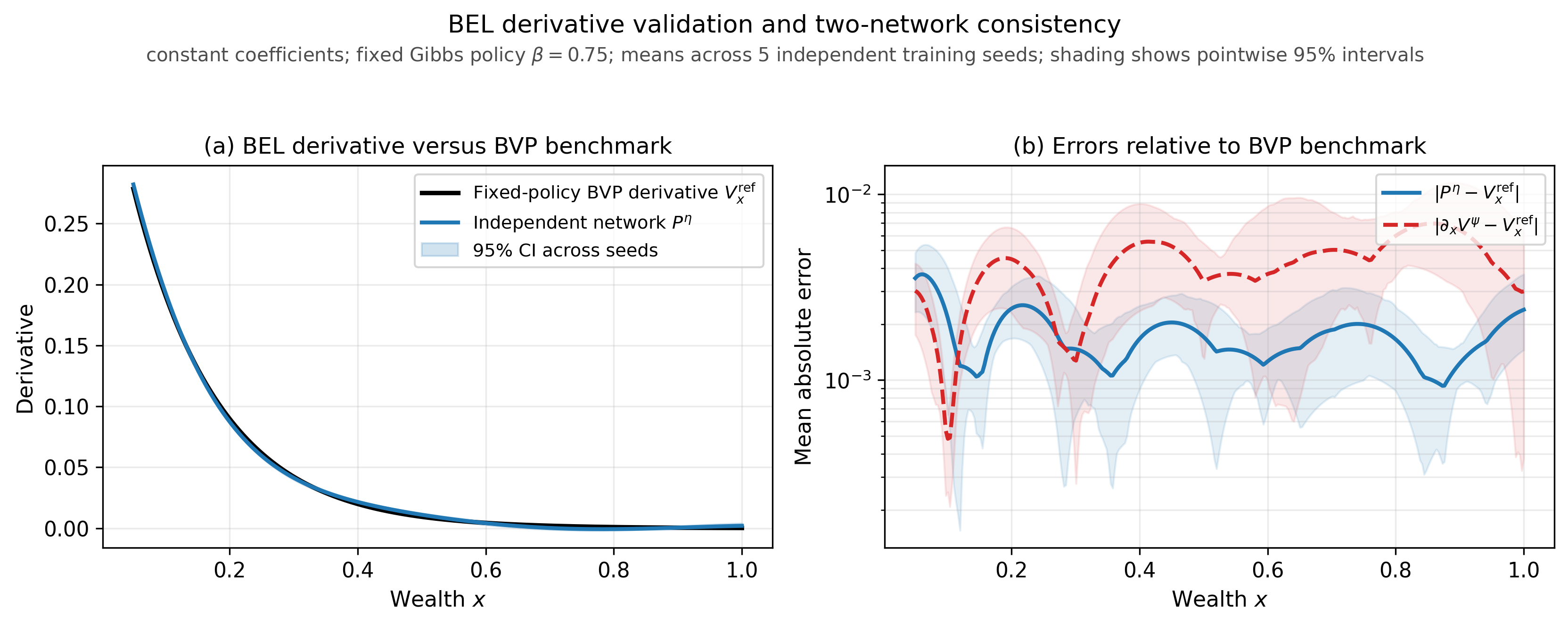}
    \caption{Validation of the independent BEL derivative approximation under
    a frozen Gibbs policy in the constant-coefficient setting.
    }
    \label{fig:bel-derivative-consistency}
\end{figure}
%Figure~\ref{fig:bel-derivative-consistency}     
In Panel~(a) above, the black curve is the deterministic derivative
obtained from the $V^{\mathrm{ref}}$, while the blue curve is the mean output of
the independently trained derivative network $P^\eta$ via the Bismut-Elworthy-Li formula.   The plot demonstrates that $P^{\eta}$ accurately represents the benchmark derivative over
$x\in[0.05,1]$.  
% Their close overlap over
% $x\in[0.05,1]$ shows that $P^{\eta}$ reproduces both the magnitude and the decay
% of the benchmark derivative. 

In Panel~(b) we evaluate the two derivative estimators separately against the benchmark $V^{\mathrm{ref}}$. The solid blue curve is
the average error $|  P^{\eta} -V_x^{\mathrm{ref}}  |$ over five seeds, 
whereas the dashed red curve is that of 
$ | \partial_xV_s^{\psi} -V_x^{\mathrm{ref}} |$ (autodiff error).
Moreover, the shaded regions represent the pointwise percentile intervals obtained
by bootstrapping the five seed-level absolute errors. The experiment shows that the $P^\eta$-$V^{\mathrm{ref}}$ error remains small throughout the wealth region $x\in[0.05,1]$ and is below the
autodiff error at essentially all grid points. Thus, the independently learned
derivative is not merely consistent with the value-network derivative; it is also the more accurate of the two estimators relative to the
external benchmark. 

More precisely, let us define the {\it Rooted Mean Square Error} between two functions $f$ and $g$ (per seed) by 
$$\mbox{RMSE}_s(f-g) := \Big[\frac{1}{N_{\mathcal G}}\sum_{x_j\in\mathcal G}\left|f_s(x_j)-g_s(x_j)\right|^2
    \Big]^{1/2}, \q s=1,\cds, 5.$$
 where $\mathcal G$ is the diagnostic grid. 
 The experiment gives the five-seed mean of RMSE as follows
\beaa
    \mbox{mean RMSE}\,
    \bigl(P^{\eta}-V_x^{\mathrm{ref}}\bigr)
    &=&1.458\times10^{-3},\\
    \mbox{mean RMSE}\,
    \bigl(\partial_xV^{\psi}-V_x^{\mathrm{ref}}\bigr)
    &=&6.553\times10^{-3},\\
    \mbox{mean RMSE}\,
    \bigl(P^{\eta}-\partial_xV^{\psi}\bigr)
    &=&6.393\times10^{-3}.
\eeaa
% Furthermore, the maximum pointwise discrepancy between the  five-seed mean  of 
% $P^{\eta}$ and $\partial_xV^{\psi}$ is $1.104\times10^{-2}$, with blue line ($|P^\eta-V_x^{\mathrm{ref}}|$-error) being consistently lower. 
This comparison numerically validates the BEL
representation and its quadrature evaluation in the present benchmark.

%We supplement Example~6.3 with two 
To end this example, we provide two sensitivity checks regarding the truncation horizon and the neural network structure. First, we vary the truncation horizon $T\in\{0.2,0.5,1,2\}$ and compare the value targets with the exact infinite-horizon fixed-policy solution. Second, with $T=2$, we use three hidden layers of widths $32$, $64$, and $128$, holding the remaining training settings fixed. At each width, the value and derivative networks are trained independently under five seeds. 
%Here the training targets and grids are fixed.) Variability across runs therefore reflects neural initialization and optimization. 
We report means and sample standard deviations across the runs.

\begin{figure}[H]
    \centering
    \includegraphics[width=0.85\textwidth]{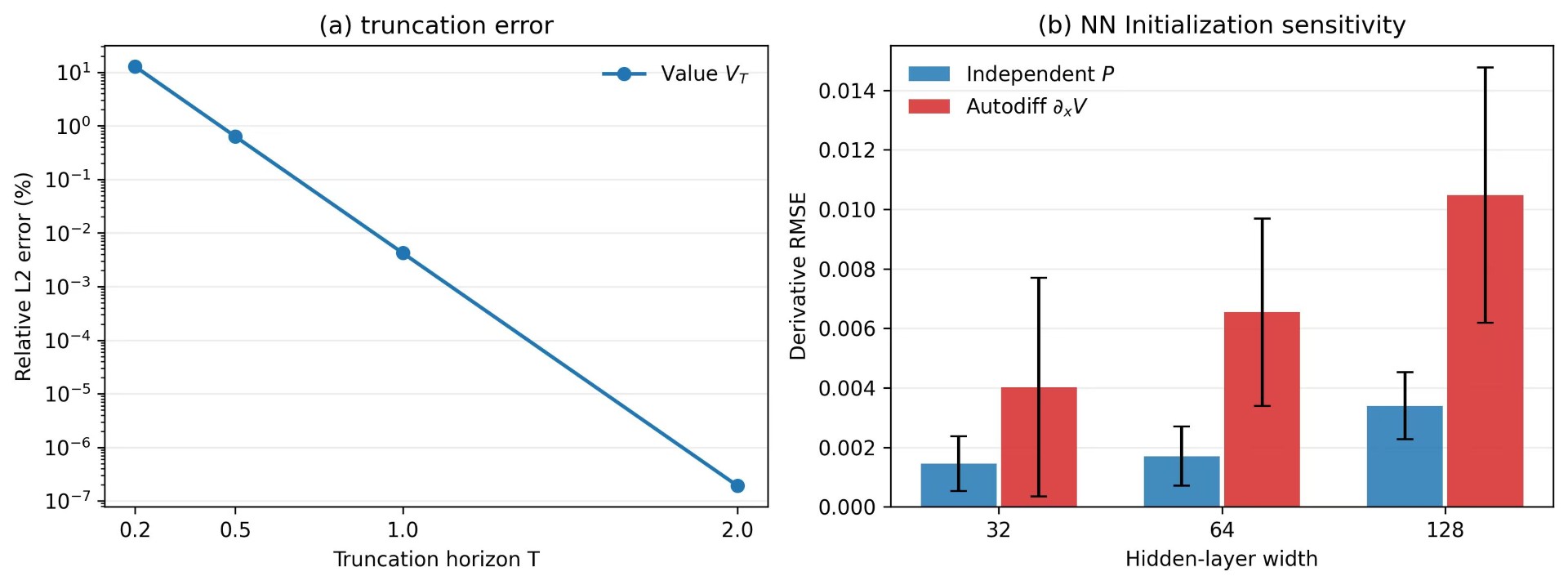}
    \caption{Sensitivity Analysis}
    \label{fig:fixed_policy_sensitivity}
\end{figure}

Figure~\ref{fig:fixed_policy_sensitivity} shows that the relative $L^2$ error of the value target (shown on a log scale) decreases from $12.83\%$ at $T=0.2$ to $0.004235\%$ at $T=1$, and is below $2\times10^{-7}\%$ at $T=2$. On the other hand, across widths $32$, $64$, and $128$, the mean RMSEs of $P^\eta$ (with benchmark derivative) are respectively $1.447\times10^{-3}$, $1.700\times10^{-3}$, and $3.395\times10^{-3}$; the corresponding RMSEs of $\partial_xV^\psi$ (with benchmark derivative) are $4.020\times10^{-3}$, $6.540\times10^{-3}$, and $1.048\times10^{-2}$. Increasing width does not improve accuracy under these optimizer settings, while the independently trained derivative network has a smaller mean error at each tested width. These observations provide insight into the robustness of the policy-iteration algorithm.}
\end{eg}

% Finally, we emphasize that the results in all the examples do not constitute a general convergence proof. They do, however,
% provide numerical evidence that the independently trained, representation-based derivative network performs better than Autodiff,  as
% an accurate approximation of the value function does not itself guarantee accurate approximation of its derivative, which is 
% %, whereas the network is trained specifically at the derivative level.does not introduce a material inconsistency. More
% %importantly, they confirm
%  the original motivation for this algorithm.
% {\color{red} Example \ref{EX3} empirically validates that $P^{\eta}$ is more accurate than
% the Autodiff derivative $\pa_xV^{\psi}$, measured against a benchmark derivative, by a
% factor of about four and a half in the five-seed mean RMSE, and that this holds
% persistently across the truncation horizons and the network widths tested. This justifies
% our scheme of learning $v'$ via an independent neural network. Examples \ref{EX1} and \ref{EX2} demonstrate
% that our PIA algorithm can be applied to approximate the optimal dividend rate with
% satisfactory empirical convergence, without knowing the model parameters, and even when
% the drift is state-dependent.
% }
{\color{black}
To summarize, Examples~\ref{EX1} and~\ref{EX2} further illustrate the application
of the proposed policy-iteration algorithm to dividend control
with constant and state-dependent drift, showing decreasing
policy-evaluation losses and convergence of the computed values empirically. Example~\ref{EX3} shows that the independently trained network
$P^{\eta}$ approximates the benchmark derivative more accurately
than the Autodiff derivative $\pa_xV^{\psi}$.
%, reducing thefive-seed mean RMSE by a factor of approximately $4.5$.
This advantage is also observed across the network widths tested,
while the truncation study supports the chosen simulation horizon
in this fixed-policy example. These findings support our approach
of learning $v'$ through a separate neural network. Together, these experiments provide numerical evidence for the
feasibility of the proposed scheme under the reported settings.
}
{\color{black}

\bibliographystyle{abbrv}
\bibliography{references}
%\bibliography{ref}

\section{Appendix}

\subsection{Proof of Proposition \ref{classic}}
Following the discussion proceeding the proposition, we know that the classical solution $v$ to (\ref{ode}) satisfying (i) exists. We need only check  (ii). 
 
 To this end, we shall follow 
%We now argue that $v$ is concave, following
 an argument of \cite{shreve}. Let us first formally  differentiate \eqref{ode}  to get
%\begin{eqnarray*}
$v'''(x)=\frac{2c}{\sigma^2}v'(x)-\frac{2}{\sigma^2}f'(v'(x))v''(x)$, $x\in[0,\infty)$.
%\end{eqnarray*}
Since  $v\in\hC^2_b([0,\infty))$, denoting $m(x):=v'(x)$, we can write 
$m''(x)=\frac{2c}{\sigma^2}m(x)-\frac{2}{\sigma^2}f'(m(x))m'(x)$, $x\in[0,\infty)$.
Now, noting Proposition \ref{fdxz}, we define a change of variables such that for $x\in[0,\infty)$,
$\varphi(x):=\int_{0}^{x}\exp\left[\int_{0}^{v}-\frac{2}{\sigma^2}f'(m(w))dw\right]dv$, 
and denote $l(y)=m(\varphi^{-1}(y))$, $y\in(0, \infty)$. Since $\varphi(0)=0$, and $\varphi'(0)=1$, we can define $\varphi^{-1}(0)=0$ as well. Then we see that,
\bea
\label{ly}
\q  \q \q l''(y)=[\varphi'(\varphi^{-1}(y))]^{-2}\frac{2c}{\sigma^2} l(y),~ y\in(0,\infty); \q l(0)=m(0)=v'(0)=  \tilde \a  >0.
\eea
Since (\ref{ly}) is a  homogeneous ODE, by uniqueness  $l(0)=  \tilde \a >0$ implies that $l(y)>0$, $y\ge 0$. That is,  $m(x)=v'(x)> 0$,  $x\ge 0$, and $v$ is (strictly) increasing. 

Finally, from  (\ref{ly}) we see that $l(y)>0$, $y\in[0,\infty)$ also implies that
%thanks to Proposition \ref{assl}, it follows  , hence 
$l''(y)> 0$, $y\in[0,\infty)$. Thus $l(\cd)$ is  convex   on $[0,+\infty)$, and hence would be unbounded unless $l'(y)\le 0$ for all $y\in[0,\infty)$. 
This, together with the fact that $v(x)$ is a bounded and increasing function, shows that $l(\cd)$ (i.e., $v'(\cd)$) 
can only be decreasing and convex, thus $v''(x)$ (i.e., 
%, as Fig.\ref{yjddtx}-(c) shows.
$l'(y)$) $\leq 0$, proving the concavity of $v$, whence the proposition.
 \qed

\subsection{Proof of Theorem \ref{conv}}

{\it  Proof.} Throughout this proof, we shall consider only the case $t\le\tilde\t_x$. Let us define, for $t\in[0,T]$ and $\th\in\Th$, 
\bea
\label{mdef}
\left\{\ba{lll}
\dis m(t,\th):=e^{-ct}J^\th(X_t)-\int_t^{\tilde\t_x}\neg\neg e^{-cs}r(X_s)\,ds, \\
\dis \hat m(t_i,\th):=e^{-ct_i}J^\th(X_{t_i})-\sum_{j=i}^{N-1}e^{-ct_j}r(X_{t_j})\D t ,
\ea\right.
\eea
for $0\le i\le N-1$, where $N=\frac{\lfloor\tilde\t_x\rfloor}{\D t}\wedge n$. By Proposition \ref{assumAP} and Assumption  \ref{ape1}, we assume that $|r(\cd)|\le R$, 
% $\ppi\in\sA^s_{cl}$, the coefficients $b^\ppi$ and
%$r^\ppi$ are bounded and Lipschitz continuous. We set $R:=\|r\|_\infty$ 
%and let $L_r$ denote the Lipschitz constant of $r$. Assumption \ref{ape1}(iii)
%provides 
$|J^\th(\cd)|\le K(\th)$, and 
% together with the Lipschitz bound
$|J^\th(y)-J^\th(z)|\le K(\th)|y-z|$, for some   $R>0$ and continuous function $K(\cd)>0$.
Therefore, it is easy to check that, for 
%. Since $t\mapsto e^{-ct}$ is decreasing,
%$\dis e^{-ct_j}\D t\le\int_{t_{j-1}}^{t_j}e^{-ct}dt$ for $j\ge1$, so that for
$\D t\le1$ and $0\le i<N$, it holds that 
%\bea
%\label{sumbound}
%\sum_{j=i}^{N-1}e^{-ct_j}\D t\;\le\;e^{-ct_i}\D t+\int_{t_i}^\infty\neg e^{-ct}dt
%\;\le\;\D t+\frac1c\;\le\;1+\frac1c .
%\eea
%Hence, since $|J^\th(\cd)|\le K(\th)$, $|r(\cd)|\le R$ and $c>0$, we have 
\bea
\label{Lam1}
\big|m(t,\th)\big| \, , \big|\hat m(t_i,\th)\big|
\;\le\;K(\th)+R\Big(1+\frac1c\Big)=:\Lambda(\th),
\qq t\in[0,\tilde\t_x],\;0\le i<N.
\eea
Next, recall \reff{ML} and  \eqref{MLtildeDt}, we can write 
\bea
\label{MLT}
\left\{\ba{lll}
\dis 2{\it ML}_T(\th)
=\hE\Big[\int_0^{\tilde\t_x}|m(t,\th)|^2dt\Big]
%\;=\;\hE\Big[\int_0^{t_N}\neg\neg\big(m(t,\th)\big)^2dt
%+\1_{\{N<n\}}\int_{t_N}^{\tilde\t_x}\neg\neg\big(m(t,\th)\big)^2dt\Big]
%\nonumber\\
=\hE\Big[\sum_{i=0}^{N-1}\int_{t_i}^{t_{i+1}}\neg\neg|m(t,\th)|^2dt
+\1_{\{N<n\}}\int_{t_N}^{\tilde\t_x}\neg\neg|m(t,\th)|^2dt\Big]; \\
%\label{MLTDt}
\dis 2 {\it ML}_{T,\D t}(\th)
=\hE\Big[\sum_{i=0}^{N-1}|\hat m(t_i,\th)|^2\D t\Big]
=\hE\Big[\sum_{i=0}^{N-1}\int_{t_i}^{t_{i+1}}\neg\neg|\hat m(t_i,\th)|^2dt\Big].
\ea\right.
\eea
Then, we have 
\bea
\label{split2}
\big|2{\it ML}_T(\th)-2{{\it ML}}_{T,\D t}(\th)\big|\;\le\;E_1(\th,\D t)+E_2(\th,\D t),
\eea
where 
\bea
\left\{\ba{lll}
\dis E_1(\th,\D t):=\hE\Big[\sum_{i=0}^{N-1}\int_{t_i}^{t_{i+1}}
\Big| \big|m(t,\th)|^2-|\hat m(t_i,\th)\big|^2\Big|dt\Big],\\
\dis E_2(\th,\D t)=\hE\Big[\int_{t_N}^{\tilde\t_x}\neg\neg
\big|m(t,\th)\big|^2\1_{\{N<n\}}dt\Big].
\ea\right.
\eea
Clearly, by \eqref{Lam1} and the fact that $\tilde\t_x-t_N\le\D t$, we have
\bea
\label{E2bound}
E_2(\th,\D t)\;\le\;\big[\Lambda(\th)\big]^2\,\D t .
\eea
To estimate $E_1(\th, t)$, we again use \eqref{Lam1} to get 
\bea
\label{E1split}
E_1(\th,\D t)
%&\le&\hE\Big[\sum_{i=0}^{N-1}\int_{t_i}^{t_{i+1}}
%\big|m(t,\th)+\hat m(t_i,\th)\big|\,\big|m(t,\th)-\hat m(t_i,\th)\big|dt\Big]
%\nonumber\\
\le 2\Lambda(\th)\sum_{i=0}^{n-1}\int_{t_i}^{t_{i+1}}
\hE\big|m(t,\th)-\hat m(t_i,\th)\big|\,dt.
\eea
Next, for each $0\le i<N$ and $t\in[t_i,t_{i+1})$, by Assumption \ref{ape1}-(iii) we can easily check that 
%Then
%\bea
%\label{tri}
%\big|m(t,\th)-\hat m(t_i,\th)\big|
%\;\le\;\big|m(t,\th)-m(t_i,\th)\big|+\big|m(t_i,\th)-\hat m(t_i,\th)\big| .
%\eea
%We have 
\bea
\label{firstterm}
|m(t,\th)-m(t_i,\th)|
&\le&=\big|e^{-ct}-e^{-ct_i}\big|\big|J^\th(X_t)\big|
+e^{-ct_i}\big|J^\th(X_t)-J^\th(X_{t_i})\big|
+\int_{t_i}^{t}\neg e^{-cs}r(X_s)\,ds \nonumber\\
%\eeaa
%so that, , and using $|e^{-ct}-e^{-ct_i}| \leq c |t-t_i|$ we have, 
%\bea
%\big|m(t,\th)-m(t_i,\th)\big|
&\le& cK(\th)|t-t_i|+K(\th)\big|X_t-X_{t_i}\big|+R|t-t_i|  \\
&\le&  L(\theta) \D t +  K(\th)\big|X_t-X_{t_i}\big|, \nonumber
\eea
where $L(\theta)= CK(\theta)+R$. Since $b^{\ppi}$ is bounded and $\sigma>0$ is constant, there exists a constant $C_0>0$ depending on $b^{\ppi}$ and $ \sigma$, such that $ \hE |X_t-X_s| \leq C_0 |t-s|^{\frac{1}{2}}$ for all $t,s$. Thus we have,  
\bea
\hE \big|m(t,\th)-m(t_i,\th)\big|
&\le&  L(\theta) \D t +  K(\th) C_0 |t-t_i|^{\frac{1}{2}}  \nonumber \\
&\le&  (L(\theta) + C_0 K (\theta)) \sqrt{\D t} \leq  \Lambda_1 (\theta)  \sqrt{\D t}
\eea
where $ \Lambda_1 (\theta):=L(\theta) + C_0 K (\theta)$.

On the other hand, we note that  
\bea
\label{secondterm}
\big|m(t_i,\th)-\hat m(t_i,\th)\big|
&=&\Big|\int_{t_i}^{\tilde\t_x}\neg\neg e^{-cs}r(X_s)ds
-\sum_{j=i}^{N-1}e^{-ct_j}r(X_{t_j})\D t\Big|
\nonumber\\
&\le&\sum_{j=i}^{N-1}\int_{t_j}^{t_{j+1}}
\big|e^{-cs}r(X_s)-e^{-ct_j}r(X_{t_j})\big|ds
+\Big|\int_{t_N}^{\tilde\t_x}\neg\neg e^{-cs}r(X_s)\1_{\{N<n\}}ds\Big|
\nonumber\\
&\le&\sum_{j=i}^{N-1}\int_{t_j}^{t_{j+1}}
\big|e^{-cs}r(X_s)-e^{-ct_j}r(X_{t_j})\big|ds+R\,\D t.
\eea
Recall that $r$ is bounded by $R>0$. Denoting $L_r$ to be the Lipschitz constant of $r$, one can easily check that 
 for   $0 \leq j \leq N-1$, and $s \in [t_j,t_{j+1}]$ one has
\beaa
\big|e^{-cs}r(X_s)-e^{-ct_j}r(X_{t_j})\big|
%&=& \big|(e^{-cs}-e^{-ct_j})r(X_{t_j})\big|+\big|e^{-cs}[r(X_s)-r(X_{t_j})]\big| \\
\le cR \D t + L_r|X_s-X_{t_j}|. 
\eeaa
Consequently, we have
%\beaa
%\hE \big|e^{-cs}r(X_s)-e^{-ct_j}r(X_{t_j})\big|
%&\leq& cR \D t + L_r B \sqrt{\D t}.
%\eeaa
%
%Thus we have 
\beaa
\label{secondterm}
\hE \big|m(t_i,\th)-\hat m(t_i,\th)\big|
&\le&\sum_{j=i}^{N-1}\int_{t_j}^{t_{j+1}}
\big(c R \D t + L_r C_0\sqrt{\D t} \big)ds+R \D t =
%\\&\le& T \big( cR \D t + L_r B \sqrt{\D t} \big)+ R \D t \leq 
\Lambda_2  \sqrt{\D t},
\eeaa
where $\Lambda_2 :=  T( cR + L_r C_0  )+ R $. It then follows that
%Thus we have 
\bea
\label{E1bound}
E_1(\th,\D t) 
&\le&2\Lambda(\th)\sum_{i=0}^{n-1}\int_{t_i}^{t_{i+1}}
\hE\big|m(t,\th)-\hat m(t_i,\th)\big|\,dt  \nonumber \\
&\le&2\Lambda(\th)\sum_{i=0}^{n-1}\int_{t_i}^{t_{i+1}}
\big( \Lambda_1(\th) \sqrt{\D t} + \Lambda_2  \sqrt{\D t}\big)\,dt = \Lambda_3(\th) \sqrt{\D t},
\eea
where $\Lambda_3(\th):= 2 \Lambda(\th) T \big(\Lambda_1(\th) + \Lambda_2\big) $.

Combining  \eqref{E2bound} and \eqref{E1bound}, and using $\Gamma(\th) = \Lambda_3(\th) + [\Lambda(\th)]^2$, we obtain the bound, \beaa
\big|{{\it ML}}_{T,\D t}(\th)-{\it ML}_T(\th)\big| \leq \Gamma(\th) \sqrt{\D t}.
\eeaa
Since $\Gamma(\th)$ is continous in $\th$, $\lim_{\D t \to 0} {{\it ML}}_{T,\D t}(\th) \to {\it ML}_T(\th) $ uniformly in $\th$ on compacta.
\qed

\end{document}